\documentclass[11pt]{article}
\usepackage{graphicx}
\usepackage{amssymb}

\newcommand{\Real}{\mathbb R}

\newcommand{\Hyp}{\mathbb{H}}
\newcommand{\Sph}{\mathbb{S}}
\newcommand{\Z}{\mathbb{Z}}
\newcommand{\Cplex}{\mathbb{C}}

\def\a{\alpha}
\def\b{\beta}

\def\e{\epsilon}

\def\C{\Gamma}

\def\eproof{$\Box$ \medskip}
\newtheorem{theorem}{Theorem}
\newtheorem{corollary}[theorem]{Corollary}
\newtheorem{lemma}[theorem]{Lemma}

\begin{document}

\title{Hausdorff dimension and the Weil-Petersson extension to quasifuchsian space\thanks{Research supported in part by NSF grant DMS 0305634}}
\author{Martin Bridgeman \\
 \footnotesize{Department of Mathematics,  Boston College,
Chestnut Hill, MA 02467}}
\maketitle
\begin{abstract}
We consider a natural non-negative two-form $G$ on quasifuchsian space that extends the Weil-Petersson metric on Teichm\"uller space. We  describe completely the positive definite locus of $G$, showing  that it is a positive definite metric off the fuchsian diagonal of quasifuchsian space and is only zero on the ``pure-bending'' tangent vectors to the fuchsian diagonal . We show that $G$ is equal to the pullback of the pressure metric from dynamics. We use the properties of $G$ to prove that at any critical point of the Hausdorff dimension function on quasifuchsian space the Hessian of the Hausdorff dimension function must be positive definite on at least a half-dimensional subspace of the tangent space. In particular this implies that Hausdorff dimension has no local maxima on quasifuchsian space.
\end{abstract}


\section{Statement of results}
\label{introsection}
Let  $S$  be a closed hyperbolic surface and $T(S)$ be the associated {\em Teichm\"uller space}. Then the {\em Weil-Petersson metric}  $w$ is a  Riemannian metric on $T(S)$.  For simplicity, we normalize the Weil-Petersson metric to define the {\em normalized Weil-Petersson metric}  $$g =  \left(\frac{2}{3\pi |\chi(S)|} \right)w\ .$$ 

If $QF(S)$ is the quasifuchsian space of $S$, then by Bers simultaneous uniformization,  $QF(S) \simeq T(S)\times T(S)$. This gives the natural diagonal embedding $\Delta:T(S) \rightarrow T(S) \times T(S) \simeq QF(S)$ given by $\Delta(X) = (X,X)$. We let $F(S) = \Delta(T(S))$ the diagonal in $QF(S)$. Then  $F(S)$ corresponds to the subspace of fuchsian elements of $QF(S)$ and is called the {\em fuchsian subspace}. It is a smooth submanifold of $QF(S)$ and we have the natural identification $T(S) \simeq F(S)$  via $\Delta$.

\medskip

Quasifuchsian  space $QF(S)$ has  complex structure coming from the fact that $QF(S)$ is an open complex submanifold of the complex representation space $$R(S) = Hom(\pi_{1}(S),PSL(2,\Cplex))/\sim$$ where $\sim$ is equivalence up to conjugation  (see \cite{Mar74}). This complex structure is  given by a bundle map $J:T(QF(S)) \rightarrow T(QF(S))$ with $J$ a lift of the identity map on $QF(S)$  and $J^{2} = -I$ where $I$ the identity map on $T(QF(S))$. 

If $h:QF(S) \rightarrow \Real$ is the map given by letting $h(X)$ be the Hausdorff dimension of the limit set of $X$, then by Ruelle \cite{Rue82}, $h$ is real-analytic.  Also associated with each $X \in QF(S)$ is a real-analytic function $L_{\mu_{X}}:QF(S) \rightarrow \Real$ by taking  the length function associated to  the {\em unit Patterson-Sullivan geodesic current} $\mu_{X}$ of $X$. 

In \cite{BT08}, we showed that the function $(h.L_{\mu_{X}})$ on $QF(S)$ is minimum  at $X$. Using this we defined a  natural  non-negative two-form $G$ on $QF(S)$ given by taking the Hessian of $(h.L_{\mu_{X}})$ at $X$. Thus
$$G_{X} =(h.L_{\mu_{X}})''(X).$$

We showed that $G$ extends the normalized Weil-Petersson metric $g$ on $F(S)$. Specifically,

\begin{theorem}{(Bridgeman-Taylor, \cite{BT08})}
There exists a continuous non-negative two-form $G$ on $QF(S)$ such that for all $X \in F(S) \subseteq QF(S)$
$$<v,w>_{G} = <v,w>_{g} \mbox{ for all } v,w \in T_{X}(F(S)) \subseteq T_{X}(QF(S)).$$
\label{naturalextension}
\end{theorem}
\vspace{-.3in}

 In this paper we answer the question of whether $G$ is a (positive-definite) metric on $QF(S)$. The answer is that $G$ is  ``almost'' a metric, in particular it is a metric off the fuchsian locus $F(S)$. The complete description of the positive-definite locus of $G$  is

\smallskip

\noindent{\bf Main Theorem}
{\em Let $v \in T_{X}(QF(S))$, $v \neq 0$.  Then $||v||_{G} = 0$ if and only if
\begin{enumerate}
\item $X \in F(S)$.
\item $v = J.w$ where $w \in T_{X}(F(S))$.
\end{enumerate}
}

Although we will not discuss this aspect further, the main theorem has a simple description in terms of the geometry of the associated deformations. If $X \in F(S) \subseteq QF(S)$ then the tangent space at $X$   decomposes into $T_{X}(QF(S)) = T_{X}(F(S)) \oplus J.T_{X}(F(S))$ (see \cite{Bon96}). If $w \in T_{X}(F(S))$ then $w$ corresponds to deforming $X$ inside the fuchsian subspace $F(S)$. Then by the earthquake theorem, this deformation is given by shearing (or twisting) $X$ along a certain measured lamination $\beta$ (see \cite{Ker85}). Thus the vectors in $T_{X}(F(S))$ are called {\em pure shearing vectors}.
If $v \in J.T_{X}(F(S))$ then  $v = J.w$ where $w$ is a pure  shearing vector with some corresponding measured lamination $\beta$. It can be shown that $v$ then corresponds to deforming the structure $X$ by bending along measured lamination $\beta$ (see \cite{Bon96}). Thus the vectors in $J.T_{X}(F(S))$ are called {\em pure bending vectors}. 

Therefore in terms of deformations, the main theorem states that the degenerate vectors for $G$ are exactly the pure bending vectors.

The proof of the main theorem is via the conformal equivalence of the two-form $G$ with another two-form $W$ obtained by taking the pullback of the so-called pressure metric  of thermodynamics. Then the proof of positive-definiteness reduces to showing  that the pullback is only trivial for the above 
tangent vectors. This  relation between  $G$ and $W$ was suggested by Curt McMullen in the paper \cite{McM06pre1}.

Using the main theorem we study properties of the critical points of $h$. In particular if $h$ is critical at $X$ then the Hessian of $h$ at $X$ is a well-defined symmetric bilinear two-form. Thus the Hessian has a well-defined signature. Applying the main theorem we obtain

\begin{theorem}
If $X \in QF(S)$ is a critical point of $h:QF(S) \rightarrow \Real$ then the Hessian of $h$ at $X$ has positive definite dimension at least $6g-6$. In particular $h$ has no local maxima in $QF(S)$.
\label{criticalpoints}
\end{theorem}
\medskip

\subsection{Background}
In \cite{BT08}, the complex structure on $QF(S)$ was used to define a ``new'' metric $H$ on Teichm\"uller space. If $X \in T_{X}(F(S))$ the associated two-form at $X$ is given by
$$<v,w>_{H} = h''(X)(J.v,J.w)$$
 for  $v,w \in T_{X}(F(S)) \subseteq T_{X}(QF(S)).$

From the definition of  $G$ and the fact that it is non-negative, we obtain
\begin{theorem}{(Bridgeman, Taylor, \cite{BT08})}
If $X \in F(S)$ and $v \in T_{X}(F(S)) \subseteq T_{X}(QF(S))$ then
$$0 \leq ||J.v||^{2}_{G} = ||v||_{H}^{2} - ||v||_{g}^{2}$$
where $g$ is the normalized Weil-Petersson metric. Thus the two-form $H$ is a positive definite metric on $F(S)$ and satisfies
$$||v||_{H} \geq ||v||_{g}.$$
\end{theorem}

In  \cite{McM06pre1}, McMullen showed that the Weil-Petersson metric was equivalent to the  the second derivative of various well-defined Hausdorff dimension  functions at the fuchsian locus. In particular McMullen proved the following theorem that the above inequality was actually an equality.
 \begin{theorem}{(McMullen, \cite{McM06pre1})}
$$ H = g.$$ 
\label{H=g}
 \end{theorem}

The results of this paper  arise out of combining the methods outlined in the paper \cite{BT08} with those of the paper of McMullen \cite{McM06pre1} and applying them in the non-fuchsian case. One note is that
we will show that $H = g$ is equivalent to the fact that $||J.v||_{G} = 0$ when $v \in T_{X}(F(S))$.

\subsection{Acknowledgements}
This paper is an outgrowth of work done in collaboration with Edward Taylor who I would like to especially thank. I would also like to especially  thank Curt McMullen for his many helpful suggestions on this project and for the paper \cite{McM06pre1}.  I would also like to thank Francis Bonahon, Dick Canary, Jeramy Kahn, and Rich Schwartz for their help.

\section{Kleinian groups and geodesic currents}
\subsection{Kleinian groups}
\hspace*{.3cm}  Let $Isom_{+}( \Hyp^{n}) \; n \geq 2$ be the space of
orientation preserving isometries of $ \Hyp^{n}.$   As is well-known, this space
of isometries can be given the topology induced by
uniform convergence on compact sets.
We define a {\em Kleinian group} $\Gamma$ to be a discrete torsion-free subgroup of $Isom_{+}( \Hyp^{n})$.  As such,
$\Gamma$ acts properly discontinuously on $ \Hyp^{n}$, and the quotient manifold
 $N =  \Hyp^{n}/\Gamma$ is a complete Riemannian manifold of constant curvature $-1$.

A Kleinian group $\Gamma$ also acts as a discrete subgroup of conformal automorphisms of the
sphere at infinity $ \Sph_{\infty}^{n-1}$; this action partitions $ \Sph_{\infty}^{n-1}$ into
two disjoint sets.  The {\em regular set} $\Omega_{\Gamma}$ is the largest open set
in $ \Sph_{\infty}^{n-1}$ on which
$\Gamma$ acts properly discontinuously, and the {\em limit set} $\Lambda_{\Gamma}$ is its complement.  In
the case that $\Lambda_{\Gamma}$ contains more than $2$ points, it is characterized as being the
smallest closed $\Gamma$-invariant subset of $ \Sph_{\infty}^{n-1}.$

Define the {\em convex hull} $CH(\Lambda_{\Gamma})$ of the limit set $\Lambda_{\Gamma}$ to be the
smallest convex subset of $ \Hyp^{n}$ so that all geodesics with both limit points in
$\Lambda_{\Gamma}$ are contained in $CH(\Lambda_{\Gamma})$.   We can take the quotient of $CH(\Lambda_{\Gamma})$
by $\Gamma$ (denoted by $C(\Gamma)$); this is the {\em convex core}.  It is the smallest
convex submanifold of $N =  \Hyp^{n}/\Gamma$ so that the inclusion map is a homotopy
equivalence.

A Kleinian group is {\em convex co-compact} if its associated convex core
is compact and it is {\em geometrically finite}
if the volume of the unit neighborhood of the convex core is finite (see Bowditch \cite{Bow93}).
This paper deals specifically with convex co-compact Kleinian groups.
For the basics in the theory of Kleinian groups we refer the reader to Maskit \cite{Mas87}.

If $\C$ is a geometrically finite Kleinian group, we define  the space $QC(\C)$ of quasiconformal deformations of $\C$ as follows;  We consider pairs $(f_{0},\C_{0})$ such that $f_{0}: \Sph^{n-1}_{\infty}
\rightarrow  \Sph^{n-1}_{\infty}$ is a quasiconformal homeomorphism,   conjugating
$\C$ to  Kleinian group $\C_{0}$, i.e. $\C_{0} = f \  \C \  f^{-1}$.  The map $f_{0}$ is called the {\em marking}. We define an equivalence relation by saying  $(f_1,\C_1) \equiv (f_2,\C_2)$  if there exists a conformal map $\a$ conjugating $\C_1$ to $\C_2$, i.e.
{
$$ f_2 \circ \gamma \circ f_2^{-1} =
(\alpha \circ f_1) \circ \gamma \circ (\alpha \circ f_1)^{-1} \ \ \mbox{ for all } \gamma \in \C.
$$
Then $QC(\C)$ is the set of equivalence classes under this equivalence relation. For convenience, we will  often supress the map $f_{0}$ in describing a point
of $QC(\C)$ and just refer to it by the group.

 \subsection{Geodesic currents}
  We can identify a  geodesic with its endpoints
on $ \Sph^{n-1}_{\infty}$ and therefore we identify the space of geodesics
on $ \Hyp^n$ by $G( \Hyp^n) \cong ( \Sph^{n-1}_{\infty} \times  \Sph^{n-1}_{\infty} -
\mbox{ diagonal})/\Z_{2}$.

 If $N$ is a convex co-compact hyperbolic $n$-manifold,  with $N =  \Hyp^{n} /\C$, then  each non-trivial homotopy class of closed curve  corresponds to a unique multiple of a primitive closed  geodesic in $N$. If $\a$ is a primitive closed geodesic in $N$, we lift $\a$ to  get a discrete subset of $G( \Hyp^n)$ which is $\C$ invariant. In this way we identify  every non-trivial homotopy class of closed curves  on $ \Hyp^n/\C$ with a  $\C$ invariant discrete subset of $G( \Hyp^n)$ and a certain integral multiplicity. We then obtain a $\C$ invariant  measure on $G( \Hyp^{n})$ by taking the Dirac measure on this discrete set times the multiplicity. This measure is the {\em geodesic current} associated with the closed curve. We have the following generalization;
\medskip

{\bf Definition:} A {\em geodesic current} for Kleinian group $\C$ is a positive measure on $G( \Hyp^n)$
that is invariant under the action of $\C$ and supported on the set of geodesics with
endpoints belonging to limit set $\Lambda_{\C}$.

As geodesic currents are (Borel $\sigma$-finite) measures, we can
add two geodesic currents and also multiply a geodesic current by
a positive constant. A geodesic current which is a constant
multiple of a closed geodesic is called a {\em discrete geodesic
current}. 

If $\C$ is a Kleinian group, we let ${\cal C}(\C)$ be
the space of geodesic currents defined  for $\C$. The natural topology on ${\cal C}(\C)$, via the
Radon-Riesz Representation Theorem, is the weak*-topology on the
space of continuous functions with compact support in $G( \Hyp^{n})$.

Below is a basic fact we will need concerning the topology on  ${\cal
C}(\C)$.  The proof involves first showing that  the geodesic flow on  the unit tangent bundle has the specification property (\cite{Bow72} and \cite{Sig74}), and then applying Theorem 1 in \cite{Sig74}.

\begin{theorem}
Let $\Gamma$ be a convex co-compact  Kleinian
group.  Then the set of discrete geodesic currents  is
dense in ${\cal C}(\C)$.
\label{dense}
\end{theorem}

We note that if $[f_{0},\Gamma_{0}] \in QC(\C)$ then $f_{0}:\Lambda_{\C} \rightarrow \Lambda_{\Gamma_{0}}$ is a homeomorphism. Therefore by pushing forward measures, we obtain a continuous homeomorphism, $f_{0}:{\cal C}(\C) \rightarrow {\cal C}(\Gamma_{0})$  (see \cite{BT05}). This map is the {\em marking} on the geodesic currents.

\subsection{Patterson-Sullivan geodesic current}
Fix $s \in  \Real^{+}.$
We define the {\em Poincar\'{e} series} of a Kleinian group $\C$ by
$$g_s(x,y) = \sum_{\gamma \in \C} e^{-s d(x,\gamma y)}$$
where $x,y \in  \Hyp^n$ and $d$ is the hyperbolic metric on $ \Hyp^n$.
Let
$$\delta(\Gamma) = \inf\{s: \;\; g_{s} < \infty\};$$
then $\delta(\Gamma)$ is called the {\em exponent of convergence} of the
Poincar\'{e} series. We refer the reader to \cite{Nic89} for further details on the exponent of convergence.

Following the work of Patterson and Sullivan, a
measure can be constructed on $ \Sph^{n-1}_{\infty}$ which is supported on
$\Lambda_{\C}$. For $x,y \in  \Hyp^n$ and $s > \delta(\Gamma)$,
we define a measure $\sigma_{x,s}$ supported on the orbit of $y$  by
$$\sigma_{x,s} =  \frac{1}{g_{s}(y,y)}\sum_{\gamma \in \Gamma} e^{-s d(x,\gamma y)}D(\gamma.y)$$ 
where $D(p)$ is Dirac measure at $p$.
The {\em Patterson-Sullivan} measure $\sigma_x$ is constructed by
taking a limit of these measures as $s \rightarrow \delta(\Gamma)^{+}$.
The measure  $\sigma_{x}$ can be
used to define a measure  $\tilde{m}$ on $(\Sph^{n-1}_{\infty}\times\Sph^{n-1}_{\infty} - \mbox{diagonal})$  given by
\begin{equation}
 d\tilde{m} = \frac{d \sigma_x (a) d\sigma_x (b)}{|b - a|^{2\delta(\Gamma)}}.
\label{psdiff}
\end{equation}
We then obtain a geodesic current $m$ by taking the pushforward of $\tilde{m}$ under  the $\Z_{2}$ cover $\pi:(\Sph^{n-1}_{\infty}\times\Sph^{n-1}_{\infty} - \mbox{diagonal}) \rightarrow G(\Hyp^{n})$ given by
$\pi(a,b) = g$ where $g$ is the geodesics with endpoints $a,b$.
This measure $m = \pi_{*}(\tilde{m})$  is $\C$-invariant and supported on  $(\Lambda_{\C} \times \Lambda_{\C} - \mbox{ diagonal})/\Z_{2}$.
Therefore it is a geodesic current and is called a {\em Patterson-Sullivan geodesic current} for $\C$.
By work of Sullivan (\cite{Sul84}), for $\C$ being geometrically
finite, $m$ is unique up to scalar multiple.

\subsection{Length functions}
Given a convex co-compact Kleinian group $\C$ then associated to each element $\gamma \in \C$ is a natural length function $L_{\gamma}:QC(\C) \rightarrow \Real$ given by letting
$L_{\gamma}([f_{0},\Gamma_{0}])$  be the translation length of the element $f_{0} \circ \gamma 
\circ f_{0}^{-1} \in \Gamma_{0}$.
This function is naturally a smooth function on $QC(\C)$. Similarly, if $\mu \in {\cal C}(\C)$ is a discrete geodesic current then $\mu$ is a multiple $r$ of a  closed geodesic  $\a$. We then  choose  $\gamma \in \C$  to be a lift of the action $\a$ and  define $L_{\mu}$   by  letting $L_{\mu} = r.L_{\gamma}$.

This can be generalized for geodesic currents to obtain the following result. 

{\bf Length Function Theorem:} {\em  (Bridgeman-Taylor, \cite{BT08}) Let $\Gamma$ be a convex co-compact Kleinian group acting on $ \Hyp^{3}$. Then there is a   continuous function 
$$L: \mathcal{C}(\Gamma) \rightarrow C^{\infty}(QC(\Gamma),\Real)$$  such that $L(\mu) = L_{\mu}$ for $\mu$ a discrete geodesic current where $C^{\infty}(QC(\Gamma),\Real)$ is the space  of smooth real-valued functions on $QC(\Gamma)$ with the $C^{\infty}$-topology.}

\medskip

Given $\mu \in {\cal C}(\C)$, we define $L_{\mu}:QC(\C) \rightarrow \Real$ by
$L_{\mu}(X) = L(\mu,X)$. The function $L_{\mu}$ is the {\em length function} for $\mu$.

We note that the continuity of $L$ implies that if $\mu_{i} \rightarrow \mu$ then $L_{\mu_{i}} \rightarrow L_{\mu}$ uniformly on compacts subsets of $QC(\C)$.

\subsection{Quasifuchsian space}
Recall that a \emph{fuchsian group} $\C$ is a finitely generated Kleinian group in $Isom_{+}(\Hyp^{3})$, with limit set $\Lambda_{\C}$ equal to a geometric circle in $\Sph^{2}_{\infty}$ and whose action preserves the components of the complement of $\Lambda_{\C}$. Identifying $\Sph^{2}_{\infty}$ with the extended complex plane $\hat{\Cplex}$,  we consider $\C$ as a group of M\"obius transformations on $\hat{\Cplex}$  with limit set equal to the extended real line $\overline{\Real}$ such that $\C$ preserves each component of $\hat{\Cplex} - \overline{\Real}$. Then the hyperbolic plane $ \Hyp^{2}$ with boundary $\overline{\Real}$ is invariant under $\C$ and $S =  \Hyp^{2}/\C$ is a hyperbolic surface.

Let $\C$ be convex co-compact and fuchsian;  we call the space
$QC(\C)$  \emph{quasifuchsian space}.  The
quotient manifold $ \Hyp^3/\C$ is homeomorphic to $S \times
 \Real$, where $S$ is the closed hyperbolic surface given by $ \Hyp^2/\C$.

To emphasize that we are dealing with a special case,  $QC(\C)$ is called the {\em quasifuchsian space} of $S$ and denoted by $QF(S)$. Also we denote the space of currents ${\cal C}(\C)$ by ${\cal C}(S)$. Furthermore we will denote the fuchsian elements of $QF(S)$ by $F(S)$.  By Bers simultaneous uniformization $QF(S) \simeq T(S) \times T(S)$ where $T(S)$ is the Teichm\"uller space of $S$ and $F(S)$ corresponds to the diagonal in $T(S)\times T(S)$  (see \cite{Bers60}). Thus if $\Delta:T(S) \rightarrow T(S)\times T(S)$ is the map $\Delta(X) = (X,X)$, then $F(S) \simeq T(S)$.

In the quasifuchsian case we have the the following extension of the real length function $L$ to a complex length function ${\cal L}$.
\medskip

{\bf Complex Length Theorem :} {\em (Bridgeman-Taylor, \cite{BT08}) For each $\mu \in {\cal C}(S)$ there exists a unique holomorphic function ${\cal L}_{\mu}:QF(S) \rightarrow  \Cplex$ with real part $L_{\mu}$ and imaginary part satisfying $Im({\cal L}_{\mu}) = 0$ on $F(S)$.
Furthermore the function $${\cal L}:{\cal C}(S) \rightarrow C^{\omega}(QF(S),  \Cplex)$$ given by ${\cal L}(\mu) = {\cal L}_{\mu}$
is continuous with respect to the topology of uniform convergence (on compacta)  on the space  $C^{\omega}(QF(S), \Cplex)$ of holomorphic functions on $QF(S)$.}
\medskip

{\bf Convention:} If $f:X \rightarrow Y$ is a smooth function then we will let $f'(x)$ denote the derivative map $f'(x):T_{x}(X) \rightarrow T_{f(x)}(Y)$. To simplify, if $v \in T_{x}(X)$ we will often write
$f'(v) = \left(f'(x)\right)(v)$. Similarly if $f'(x) = 0$ then the Hessian of $f$ is denoted by $f''(x)$ and is  the well-defined symmetric bilinear two-form given by
$$\left(f''(x)\right)(v,w) = \frac{\partial^{2}{f}}{\partial v \partial w}.$$
Once again we will often shorten and write $f''(v,w) = (f''(x))(v,w)$.

\section{Weil-Petersson extension $G$}
We now describe the symmetric bilinear form $G$ on $QF(S)$ given in \cite{BT08}.

Let $X = [f_{0},\Gamma_{0}] \in QF(S)$, then $f_{0}$  gives a natural homeomorphim  $f_{0}: {\cal C}(S)= {\cal C}(\C) \rightarrow {\cal C}(\Gamma_{0})$  between geodesic current spaces coming from the marking. We let $m_{\Gamma_{0}} \in {\cal C}(\Gamma_{0})$ be a Patterson-Sullivan geodesic current and pullback to define $m_{X} = f^{-1}_{0}(m_{\Gamma_{0}}) \in {\cal C}(S)$. We normalize  to define the unit length Patterson-Sullivan geodesic current of $X$ by 
$$\mu_{X} = \frac{m_{X}}{L(X,m_{X})}.$$
Then  this geodesic current has unit length in $X$.
  
In \cite{BT08}, we show that the function $(h.L_{\mu_{X}}): QF(S) \rightarrow \Real$ given by $(hL_{\mu_{X}})(Y) = h(Y).L_{\mu_{X}}(Y)$ is minimum  at $X$. Using this  we defined $G$  to be the non-negative two-form  with symmetric bilinear form  at $X$  given by
$$G_{X}  = (h.L_{\mu_{X}})''(X).$$ Finally we proved theorem \ref{naturalextension}, showing that $G$ is a natural extension of the normalized Weil-Petersson metric on $F(S)$.

\section{Thermodynamics and pressure metric}
We will now describe the pressure metric for a shift of finite type. This will be a cursory introduction to the elements of Thermodynamic Formalism needed to state and prove our results. For a complete description see the book \cite{PP90} by Parry and Pollicott and the  paper \cite{McM06pre1} of McMullen.

Let  $A$ be  a $k\times k$ matrix of zeros and ones then we define the associated {\em (one-sided) shift of finite type}  by $(\Sigma,\sigma)$ where $\Sigma$ is the set of sequences
$$\Sigma = \left\{ x = (x_{n})_{n=0}^{\infty}: x_{n} \in \{1,\ldots,k\}, A(x_{n},x_{n+1}) = 1\right\}$$
and $\sigma:\Sigma \rightarrow \Sigma$ is the standard shift where $\sigma(x_{0},x_{1},x_{2},\ldots) = (x_{1},x_{2}, \ldots)$.   We give $\{i,\ldots,k\}$ the discrete topology and $\Sigma$ the associated product topology. 

The space $C(\Sigma)$ is the space of continuous real valued functions on $\Sigma$. Two function $f,g \in C(\Sigma)$ are {\em cohomologous} ($f \sim g$),  if there exists a continuous function $h \in C(\Sigma)$ such that $f(x) -  g(x) =  h(\sigma(x)) - h(x)$. If $f \sim 0$ then $f$ is a {\em coboundary}.

We can metrize the topology on $\Sigma$  by choosing any $ K > 1$ and then defining  $d(x,y) = K^{-N}$ where $N  = N(x,y)= \min\{n\ |\  x_{n} \neq y_{n}\}$. 

Then given  $\theta \in (0,1)$ we say $f \in F_{\theta}(\Sigma)$ if there exists a constant $C > 0$ such that
$$|f(x) - f(y)| \leq C.\theta^{N(x,y)}.$$
 The set $F_{\theta}$ is the set of H\"older continuous functions with the same H\"older constant, with respect to the metric $d$.

$F_{\theta}(\Sigma)$ is given the norm $||.||_{\theta}$ by
$$||f||_{\theta} = ||f(x)||_{\infty} + \sup_{x \neq y}\frac{|f(x)-f(y)|}{\theta^{N(x,y)}}$$

Given a map $f$  we can take the iterated sum $S_{n}f$ defined by
$$(S_{n}f)(x) = \sum_{k=0}^{n-1}f(\sigma^{k}(x)).$$
If $f \sim g$ with $f(x) - g(x) = h(\sigma(x)) - h(x)$ then
$$S_{n}f(x) - S_{n}g(x) = h(\sigma^{n}(x)) - h(x).$$
Also if $f \in  F_{\theta}(\Sigma)$, the {\em Ruelle operator} $L_{f}:F_{\theta}(\Sigma) \rightarrow F_{\theta}(\Sigma)$ is defined by
$$(L_{f}g)(x) = \sum_{\sigma(y) = x} e^{f(y)}g(y).$$
We note that under iteration of the Ruelle operator we have 
$$(L^{n}_{f}g)(x) = \sum_{\sigma^{n}(y) = x}e^{S_{n}f(y)}g(y).$$

The shift $(\Sigma,\sigma)$ is {\em aperiodic} if there exists an $n >0$ such that $A^{n}$ is all positive entries. We have the following generalization of the Perron-Frobinius theorem for matrices.
 
\begin{theorem}{(Ruelle-Perron-Frobinius, \cite{PP90})} 
Let $f \in F_{\theta}(\Sigma)$ and $(\Sigma,\sigma)$ be aperiodic. Then
\begin{enumerate}
\item There is a simple maximal positive eigenvalue $\beta$ for $L_{f}$ with corresponding strictly positive eigenvector $h$.
\item The remainder of the spectrum of $L_{f}$ is contained in a disk of radius strictly smaller than $\b$. 
\item There is a unique probability measure $\mu$ such that $L_{f}^{*}\mu = \beta  . \mu$.
\item Let $h$ be  a maximal eigenvector normalized so that   $\mu(h) = 1.$ Then 
$$ \frac{L_{f}^{n}(g)}{\b^{n}}  \rightarrow h.\int gd\mu \mbox{ 
uniformly for all }g \in C(\Sigma).$$
\end{enumerate}
\label{rpf}
\end{theorem}

The {\em pressure } $P(f)$  is defined by $P(f) = \log \beta$. If $f \in F_{\theta}(\Sigma)$ satisfies $P(f) = 0$ and $h$ is a maximal normalized eigenvector of $L_{f}$ then the measure $m  = h.\mu$ is an ergodic $\sigma$-invariant probability measure and is called 
the {\em equilibrium measure} of $f$.

In \cite{PP90} the properties of the function $P:F_{\theta}(\Sigma) \rightarrow F_{\theta}(\Sigma)$ are described in detail. In particular it is  convex and real-analytic and depends only on cohomology class. 

Also if $P(f) = 0$, with equilibrium measure $m$ and $g \in F_{\theta}(\Sigma)$ then
$$P'(f)(g) = \frac{d}{dt} P(f+tg)|_{t=0} = \int g\ dm$$

Also if $P'(f)(g) = 0$ then the {\em variance} $Var(g, m)$ is defined by
 $$P''(f)(g) = \frac{d^{2}}{dt^{2}} P(f+tg)|_{t=0} = Var(g,m).$$

We define
$$T(\Sigma) =\left\{ f: f \in F_{\theta}(\Sigma) \mbox{ some } \theta , P(f) = 0\right\}/\sim.$$
Then $T(\Sigma)$ is the set of pressure zero, H\"older continuous functions up to co-boundary. If $[f] \in T(\Sigma)$ and $f$ has equilibrium measure $m$, then by the formula for the derivative of pressure $P$, the tangent space of $T(\Sigma)$ at $[f]$ can be identified with
$$T_{[f]}T(\Sigma) = \left\{ g: \int g\  dm = 0 \right\}/\sim.$$
The pressure metric $||.||_{P}$ on $T(\Sigma)$ is then defined by
\begin{equation}
||[g]||_{P} = \frac{Var(g, m)}{-\int f\   dm}.
\label{pressuremetricdef}
\end{equation}

By  Theorem 4.2  of \cite{PP90}, $Var(g, m) = 0$ implies that $g \sim 0$. Thus $||[g]||_{P} = 0$ implies $[g] = 0$ and therefore $||.||_{P}$ is positive definite metric on $T(\Sigma)$.

\section{Thermodynamics on $QF(S)$}
Let $\Gamma$ be a Kleinian group with limit set $\Lambda_{\Gamma} \subset \hat{\Cplex}$.
A conformal {\em Markov} map for $\Gamma$ is a piecewise  conformal map $f:\Lambda_{\Gamma} \rightarrow \Lambda_{\Gamma}$ such that $\Lambda_{\Gamma}$ has a partition into segments $J_{1}, \ldots, J_{m}$ so that;
\begin{enumerate}
\item $   f|_{J_{k}}  = \gamma_{k}|_{J_{k}} \mbox{ for some } \gamma_{k} \in \Gamma$
\item for each $k$, $f(J_{k})$ is the union of various $J_{l}$'s.
\end{enumerate}
A Markov map is {\em expanding} if there is an $n > 0$ such that the n-th iterate $f^{n} = f \circ f \circ\ldots \circ f$ has derivative  whose length in the spherical metric  satisfies
$$|(f^{n})'(x)| > C > 1$$ and for any $U \subset L_{\Gamma}$ open, there exists an $m > 0$ such that $f^{m}(U) = \Lambda_{\Gamma}.$

If $\Gamma$ has an expanding Markov map $f$ then we can define a matrix $A$ by $A(i,j) = 1$ if 
$J_{j} \subset f(J_{i})$ and zero otherwise. Then we have an aperiodic shift $(\Sigma,\sigma)$ and we define $\pi:\Sigma \rightarrow L_{\Gamma}$ by $\pi(x) = z$ where $f^{i}(z) \in J_{x_{i}}$. The map $f$ obviously satisfies $f(\pi(x)) = \pi(\sigma(x))$.  The map $\pi$ is surjective but as the segments $J_{i}$ may have endpoints in common, the map $\pi$ is two to one on a countable set of points $P$. If $Q$ is the finite set of endpoints of the $J_{i}$'s then $P$ is precisely
$$P = \bigcup_{n=0}^{\infty} f^{-n}(Q).$$
The points of $P$ are called {\em bad} points and if $z \not\in P$ it is called a {\em good} point. We note hat if $z$ is a good point, then there is a unique $x \in \Sigma$ such that $\pi(x) = z$ and for any $n > 0$, there is unique $\gamma_{n} \in \Gamma$  such that $f^{n} = \gamma_{n}$ on an open interval about $z$.


\subsection{Expanding Markov map for quasifuchsian groups \label{Bowensection}}
In the following we describe Bowen's results from \cite{Bow79} on  expanding Markov maps for quasifuchsian groups. 

Bowen first considered the co-compact fuchsian group  $\C_{r}$ obtained by  identifying  sides of a regular hyperbolic $4n$-gon in the standard way given by the side labelling
$$x_{1}y_{1}x^{-1}_{1}y^{-1}_{1}\ldots x_{n}y_{n}x^{-1}_{n}y^{-1}_{n}.$$
He then described an expanding Markov map $f_{\C_{r}}:\Sph^{1} \rightarrow \Sph^{1}$ for  $\Gamma_{r}$ which we will describe in detail below. 

Then if   $g:\hat{\Cplex} \rightarrow \hat{\Cplex}$  is a quasiconformal map conjugating the action of $\C_{r}$ to the action of $\Gamma$, then this gives the map $f_{\Gamma}: \Lambda_{\Gamma} \rightarrow \Lambda_{\Gamma}$ by
$f_{\Gamma} = g \circ f_{\C_{r}} \circ g^{-1}$ and $\pi_{\Gamma}: \Sigma \rightarrow \Lambda_{\Gamma}$ by $\pi_{\Gamma} = g \circ \pi_{\C_{r}}$. Then $f_{\Gamma}$ is an expanding Markov map for $\Gamma$ with the same shift space $(\Sigma, \sigma)$.

The function $\phi_{\Gamma}:\Sigma \rightarrow \Real$  defined by $\phi_{\Gamma}(x) = -\log | f'_{\Gamma} (\pi_{\Gamma}(x))|$ is H\"older continuous. By the chain rule for differentiation we have
\begin{equation}
(S_{n}\phi_{\Gamma})(x) = - \log|(f^{n}_{\Gamma})'(\pi_{\Gamma}(x))|.
\label{iteration}
\end{equation}

Then if $h_{\Gamma}$ is the Hausdorff dimension of the limit set $\Lambda_{\Gamma}$, Bowen showed  that $h_{\Gamma}$ is characterized by the equation
\begin{equation}
P(h_{\Gamma}.\phi_{\Gamma}) = 0.
\label{Pchar}
\end{equation}

We  now describe the map $f_{\C_{r}}$ in more detail.  The group $\C_{r}$ has fundamental domain $D$, the regular hyperbolic $4n$-gon. We label the sides of $D$  by $s_{i}, i = 1,\ldots,4n$. Each $s_{i}$ belongs to a unique geodesic $g_{i}$ with endpoints $p_{i},q_{i}$ on $\Sph^{1}$ . We let $I_{i}$ be the interval on $\Sph^{1}$ with endpoints $p_{i}, q_{i}$ which is smallest in length. We further define $\gamma_{i} \in \C_{r}$ to be the element which identifies $s_{i}$ with another side $s_{j}$ of $D$ for some $j$.  

For each $\gamma \in \C_{r}$ we let $D_{\gamma} = \gamma(D)$  and say $D_{\gamma}$ {\em abuts} $D$ if $D_{\gamma} \cap D \neq \emptyset$. For each $D_{\gamma}$ we define the intervals $I_{\gamma,i} = \gamma(I_{i})$. We let $R \subseteq \Sph^{1}$ be the union of the endpoints of $I_{\gamma,i}$ for $D_{\gamma}$ abutting $D$. Then $R$ defines a decomposition of $\Sph^{1}$ into intervals $J_{k}$. We note that each $J_{k} \subseteq I_{j}$ for some $j$, not necessarily unique.
We then define 
$$f_{\C_{r}}|_{J_{k}} = \gamma_{j} \mbox{ where } J_{k} \subseteq I_{j}.$$
It follows that for any choice of $j$ such that $J_{k}\subseteq I_{j}$, the map $f_{\C_{r}}$ is a Markov map for $\C_{r}$.  Bowen makes a specific choice to define $f_{\C_{r}}: \Sph^{1} \rightarrow \Sph^{1}$ such that $f_{\C_{r}}$ is an expanding Markov map and  orbit equivalent to the action of  $\C_{r}$ on $\Sph^{1}$ (see \cite{Bow79}).

 We say a geodesic $g$ {\em abuts} $D$ if it intersects a domain that abuts $D$. We now describe an elementary property of the map $f_{\C_{r}}$ that follows easily from its definition.

\begin{lemma}
Let $f_{\C_{r}}:\Sph^{1} \rightarrow \Sph^{1}$ be the expanding map for the $\C_{r}$ described above. Let $g$ be a geodesic with endpoints $a,b$  that abuts $D$. If $f_{\C_{r}} = \gamma \in \C_{r}$ at $a$, then
geodesic $\gamma(g)$ abuts $D$.
\label{goodgeodesics}
\end{lemma}

\begin{figure}
\begin{center}
\includegraphics[width=3in]{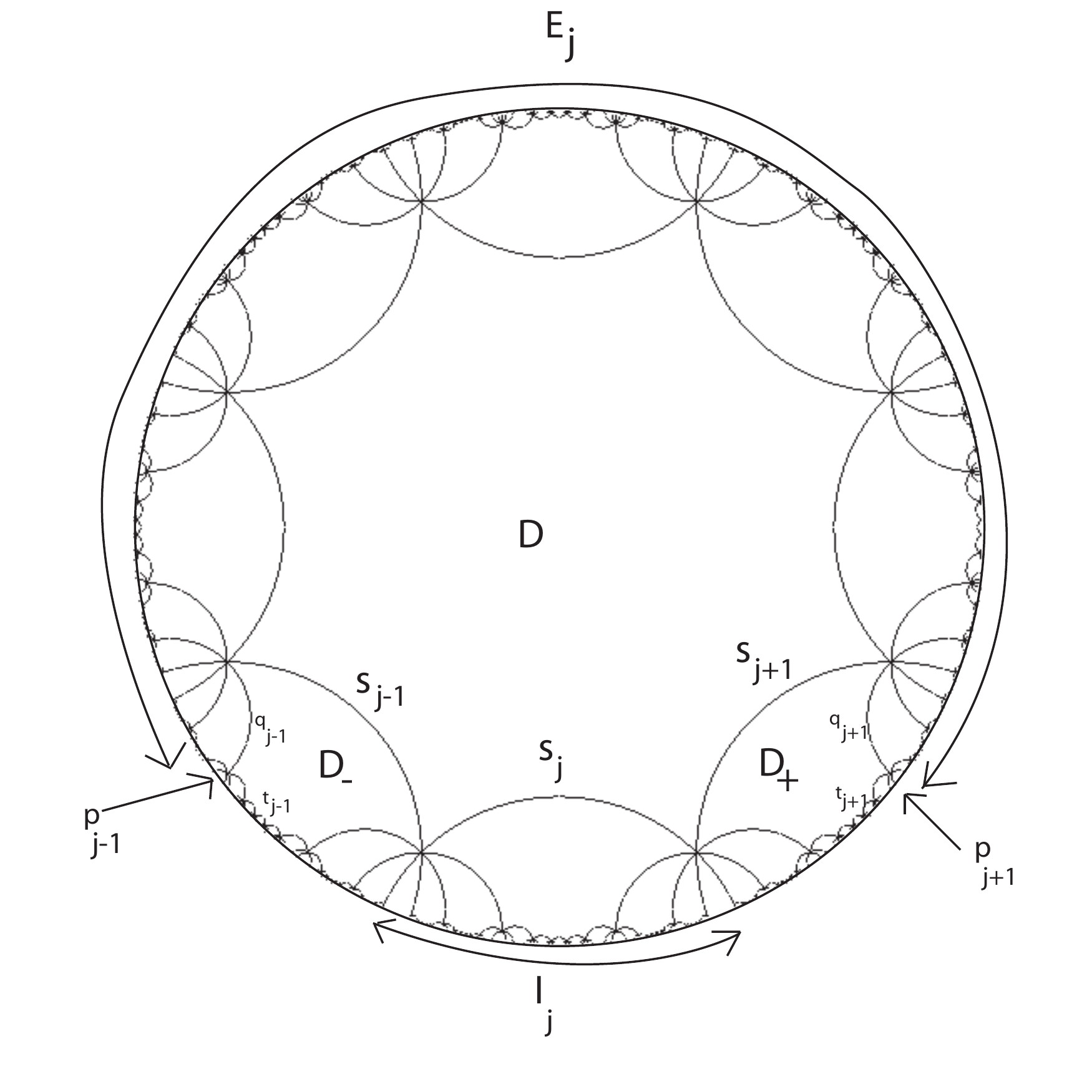}
\caption{Tesselation by regular 4n-gons}
\label{tess}
\end{center}
\end{figure}

{\bf Proof:}
We place $D$ with the origin in the center. We label the edges $e_{i}$  of $D$ clockwise for $i = 1,\ldots,4n$. The edges $e_{i}$ define geodesics $g_{i}$ which given overlapping intervals  $I_{i}$. We further define the half-plane given by $I_{i}$, $H_{i}$ (see figure 1).

Let $g$ be a geodesic with endpoints $a,b$ which abuts $D$ and $f_{\C_{r}} = \gamma $ at $a$.
Then $\gamma = \gamma_{j}$ for some $j$ and $a \in I_{j}$. We let $P$ be the convex polygon obtained by taking the union of all domains that abut $D$.  Then by assumption $g$ intersects $P$. 

If $s_{j}$ is the side of $D$ corresponding to the side identification $\gamma_{j}$, we let $P_{1}$ be the collection of domains which intersect $s_{j}$ (i.e. contain $s_{j}$ as a side or contain an endpoint of $s_{j}$ as a vertex). As $\gamma_{j}(s_{j}) = s_{k}$ some other side of $D$, if $g$ intersects $P_{1}$ then $\gamma_{j}(g)$ abuts $D$. We will now show that if $g$ intersects $P$ then it intersects $P_{1}$ thereby proving the result.

Let $P_{2}$ be the set of  domains  $P \cap H_{j}$. As $H_{j} \cap D = s_{j}$, then $P_{2 } \cap P_{1}$. As we are assuming $g$ does not intersect $P_{1}$, then $g$ does not intersect $P_{2}$. If  both endpoints of  $g$ are in $I_{j}$ then by convexity $g \subseteq H_{j}$ and therefore $\emptyset \neq g\cap P = g \cap (P\cap H_{j}) = g \cap P_{2} \subseteq g \cap P_{1}$. Thus if $g$ has both endpoints in $I_{j}$ then $g$ intersects $P_{1}$.

We let $D_{-}$ be the domain that shares the side $s_{j-1}$ with $D$ and $D_{+}$ be the domain that shares the side $s_{j+1}$ with $D$. Then as $s_{j-1}, s_{j+1}$  share a vertex with$s_{j}$ then $D_{-}, D_{+} \subseteq P_{1}$.  

In domain $D_{-}$ we label the opposite side to $s_{j-1}$ by $t_{j-1}$.  We let $h_{j-1}$ be the geodesic associated to $t_{j-1}$ and the interval $T_{j-1} \subseteq I_{j-1}$. The interval $T_{j-1}$ does not intersect with any other interval $I_{k}$ for $k \neq j-1$. Also $D_{-}$ shares a unique side with a domain 
of $P-P_{1}$ (abutting $D$ but not side $s_{j}$). We label the geodesic to the side $q_{j-1}$ and its unique endpoint $p_{j-1} \in I_{j-1}$.
We define similar quantities for $D_{+}$. We let $E_{j}$ be the interval in $\Sph^{1}$ with endpoints $p_{j-1}, p_{j+1}$ and not containing $I_{j}$. Then $E_{j}$ is disjoint from $T_{j-1}$ and $T_{j+1}$. If $g$ intersects a domain of $P$ but not of $P_{1}$ then $g$ must have an endpoint in $E_{j}$. Therefore  $g$ separates the geodesics $h_{j-1}, h_{j+1}$. Therefore $g$ must intersect either $D_{-}, D_{+}$ or $D$. Thus $g$ intersects $P_{1}$ and we're done.
\eproof

The above lemma says that the abutting geodesics are an invariant set under the M\"obius map defined by their endpoints.

We now prove  an important property of the expanding Markov map $f_{\C_{r}}$ that we will need later.

If $G$ is a group, then we say $g$ is commensurable to  $h$ if there exists $k \in G$ such that $g^{n} = k h^{m}k^{-1}$  for some $n,m \neq 0$. The set of commensurability classes of $G$ is denoted $[G]$. We note that for $\Gamma$ a co-compact Kleinian group, then $[\Gamma]$ is equivalent to the set of primitive geodesics in $N = \Hyp^{n}/\Gamma$.

\begin{lemma}
 Let  $f_{\C_{r}}:\Sph^{1} \rightarrow \Sph^{1}$ be the expanding Markov map  described  above. Then there is  a finite set $S \subseteq [\C_{r}]$ of commensurability classes  such that if  $[\gamma] \not
\in S$, then 
\begin{enumerate}
\item If $\gamma' \in [\gamma]$ then the endpoints of the axis of $\gamma'$ are good points of $f_{\C_{r}}$. 
\item There exists a $\gamma' \in [\gamma]$  whose axis abuts $D$ and has  fixed points $a, b$ such that  the expanding fixed point $a$ of $\gamma'$ is a periodic point of $f_{\C_{r}}$ \end{enumerate}
\label{techlemma}
\end{lemma}

{\bf Proof:} We note that as $\C_{r}$ does not contain parabolics,  then if two elements $\gamma_{1}, \gamma_{2}$ have a common fixed point, then they share the same axis and are commensurate. Let $D$ be the fundamental domain for $\C_{r}$ given by the regular $4n$-gon in the Poincare disk model with center at $0$. We extend the $4n$ sides of $D$ to complete geodesics $g_{i}, i = 1,\ldots,4n$, and let $P_{0}$ be the union of the endpoints of the $g_{i}$'s. We then define $P = \C_{r}.P_{0}$, the orbit of $P_{0}$ under the group. The set $P$ is precisely the set of bad points for $f_{\C_{r}}$. Also as $D$ is the regular $4n$-gon, each of the geodesics $g_{i}$ is the axis for an element of $\C$ and we let $\gamma_{i}$ be the corresponding primitive element.  Then we let $S = \{[\gamma_{i}]\}_{i=1}^{4n}$. If an element $\gamma \in \C_{r}$ has a bad endpoint $z \in P$, then $z = \gamma_{1}.z_{0}$ for $z_{0} \in P_{0}$ then $\gamma_{1}^{-1}\gamma \gamma_{1}$ has fixed point $z_{0}$. Therefore $\gamma_{1}^{-1}\gamma \gamma_{1}$ shares an endpoint with some $\gamma_{i}$. Therefore $\gamma$ and $\gamma_{i}$ are commensurate. Thus we conclude that if $[\gamma] \not\in S$, then $\gamma$ has both endpoints being good.

We now let $[\gamma] \not\in S$ and choose $\gamma$ such that its axis intersects $D$. As $[\gamma] \not\in S$, then its axis $g$ has endpoints $a_{0},b_{0}$. 
We define  geodesic $g_{n}$ to have endpoints $a_{n} = f_{\C_{r}}^{n}(a)$, and let $b_{n} \in \Sph^{1}$ be the unique point such that the pair  $(a_{n},b_{n}) \in \Sph^{1} \times \Sph^{1}$ are endpoints of the axis of a conjugate of $\gamma$.  Then by lemma \ref{goodgeodesics}, $g_{n}$ also abuts $D$. By compactness of the union of domains abutting $D$, we have there is an $\e > 0$ such that $|a_{n} - b_{n}| > \e$ for all $n$.  But if the sequence $\{(a_{n}, b_{n})\}$ in  $\Sph^{1} \times \Sph^{1}$ has an infinite number of values, it must have a convergent subsequence. As the  orbit of the axis $g$ under $\C_{r}$ is discrete in the space of geodesics $G(\Hyp^{2})$, any convergent subsequence  must converge to a point on the diagonal of $\Sph^{1}\times \Sph^{1}$ contradicting $|a_{n} - b_{n}| > \e$. Therefore the sequence takes a finite set of values and  there exists a $k > 0$ and an $n$ such that $(a_{n},b_{n}) = (a_{n+k},b_{n+k})$ and therefore $f^{k}_{\C_{r}}(a_{n}) = a_{n}$ and $a_{n}$ is a periodic point for $f_{\C_{r}}$. We let $\gamma'$ be the conjugate of $\gamma$ with endpoints $(a_{n},b_{n})$, giving the result.
\eproof

\subsection{Pullback of pressure metric}

We first define the {\em Hausdorff dimension function} $h:QF(S) \rightarrow \Real$ given by letting $h([f_{0},\Gamma_{0}]) = h_{\Gamma_{0}}$, the Hausdorff dimension of the limit set $\Lambda_{\Gamma_{0}}$. This is  well-defined, and by Ruelle (see \cite{Rue82}), $h$ is real-analytic. 

We let $\C$ to be a fuchsian group such  that $S = \Hyp^{2}/\C$ with  expanding Markov map $f_{\C}$ as described in the section \ref{Bowensection}. We let $(\Sigma,\sigma)$ be the associated shift and $\pi_{\C};\Sigma \rightarrow \Sph^{1}$ as before. Then for each $X \in QF(S)$ we let $X = [g_{0},\Gamma_{0}]$ where $g_{0}$ conjugates $\C$ to $\Gamma_{0}$. 
Then we define $\phi_{X} = \phi_{\Gamma_{0}}: \Sigma \rightarrow \Real$. 

\begin{lemma}
If $[g_{0},\Gamma_{0}] =[g_{1},\Gamma_{1}] \in QF(S)$ then $\phi_{\Gamma_{0}} \sim \phi_{\Gamma_{1}}$.
\label{Fwelldefined}
\end{lemma} 
{\bf Proof:}
We note that if $(g_{0},\Gamma_{0}) \sim (g_{1},\Gamma_{1})$ then there is a conformal map $c:\hat{\Cplex} \rightarrow :\hat{\Cplex}$ conjugating $\Gamma_{0}$ to $\Gamma_{1}$. Therefore $c \circ f_{\Gamma_{0}} = f_{\Gamma_{1}} \circ c$, giving
$$ c'(f_{\Gamma_{0}}(z))f'_{\Gamma_{0}}(z) =  f'_{\Gamma_{1}}(c(z)).c'(z).$$
As $c$ conjugates the action of $\Gamma_{0}$ to the action of $\Gamma_{1}$, we have $c\circ g_{0} = g_{1}$ on the limit set $\Lambda_{\Gamma_{0}}$. Therefore $c \circ \pi_{\Gamma_{0}} = c\circ g_{0} \circ \pi= g_{1}\circ \pi =  \pi_{\Gamma_{1}}$. Therefore if  $z = \pi_{\Gamma_{0}}(x)$ then
$$ c'(f_{\Gamma_{0}}(\pi_{\Gamma_{0}}(x)))f'_{\Gamma_{0}}(\pi_{\Gamma_{0}}(x)) =  f'_{\Gamma_{1}}(c(\pi_{\Gamma_{0}}(x))).c'(\pi_{\Gamma_{0}}(x)) =  f'_{\Gamma_{1}}(\pi_{\Gamma_{1}}(x)).c'(\pi_{\Gamma_{0}}(x)).$$
Taking logs of absolute values, we have
$$ \phi_{\Gamma_{0}}(x) - \phi_{\Gamma_{1}}(x) = \log|c'(f_{\Gamma_{0}}(\pi_{\Gamma_{0}}(x))| - \log|c'(\pi_{\Gamma_{0}}(x))|$$
By definition, we have $f_{\Gamma_{0}} \circ \pi_{\Gamma_{0}} = \pi_{\Gamma_{0}}\circ \sigma$. Therefore we let $q \in C(\Sigma)$ be given by $q(x) = \log|c'(\pi_{\Gamma_{0}}(x))|$. Then
 $$ \phi_{\Gamma_{0}}(x) - \phi_{\Gamma_{1}}(x) = q(\sigma(x)) - q(x)$$
and $\phi_{\Gamma_{0}}, \phi_{\Gamma_{1}}$ are cohomological. \eproof

We now define $\Phi_{X} = h_{\Gamma_{0}}.\phi_{\Gamma_{0}}$. Then by equation \ref{Pchar}, $P(\Phi_{X}) = 0$. We obtain a map
$F:QF(S) \rightarrow T(\Sigma)$ by $F(X) = [\Phi_{X}]$. By the above lemma \ref{Fwelldefined},  the map $F$ is well-defined. We then define $W$ to be the pullback of the pressure metric on $T(\Sigma)$. As the pressure metric is positive definite it follows that $W$ is at least non-negative. 

 To obtain a formula for $||.||_{W}$, given a $v \in T_{X}(QF(S))$, we  choose a smooth curve $\a:(-\e,\e) \rightarrow QF(S)$ with $\a(0) = X$ and $\a'(0) = v$. Then let $\a(t)= X_{t} = [g_{t},\Gamma_{t}]$ where $g_{t}$ is a smooth 1-parameter family of quasiconformal maps conjugating $\C$ to $\Gamma_{t}$.

We let $\Phi_{t} = \Phi_{\a(t)}$ and $\phi_{t} = \phi_{\a(t)}$ and define
$\dot{\Phi}_{0}$ by
$$\dot{\Phi}_{0}(x) = \frac{d}{dt}\left|_{t=0}\left(\Phi_{t}(x)\right)\right.$$
Then by definition of the pressure metric in equation \ref{pressuremetricdef}, $||v||_{W} $ is given by
$$||v||_{W}^{2} = \frac{Var(\dot{\Phi}_{0}, m)}{-\int \Phi_{0} \  dm}.$$

We obtain an alternative definition of $||.||_{W}$ by noting that $P(\Phi_{t}) = 0$ and taking  derivatives with respect to $t$. 
Taking first derivatives we obtain
$$P'(\Phi_{t})(\dot{\Phi}_{t}) = 0.$$
Then taking derivative again we have
$$ P''(\Phi_{t})(\dot{\Phi}_{t}) + P'(\Phi_{t})(\ddot{\Phi_{t}}) = 0.$$ 
Evaluating at $t = 0$ we have
$$Var(\dot{\Phi}_{0},m) + \int \ddot{\Phi}_{0}\  dm = 0.$$
Therefore we have
\begin{equation}
||v||^{2}_{W} = \frac{Var(\dot{\Phi}_{0}, m)}{-\int \Phi_{0}\ dm} =  \frac{\int \ddot{\Phi}_{0}\ dm}{\int \Phi_{0} \ dm}.
\label{metricdefalt}
\end{equation}

\section{Conformal equivalence of $G$ and $W$}
The proof that $G$ and $W$ are conformally equivalent follows by generalizing  the argument in \cite{McM06pre1} of McMullen  for the fuchsian subspace $F(S)$  to all of quasifuchsian space $QF(S)$.

\begin{theorem}
The pseudometrics $G$ and $W$ are conformally equivalent with 
 $$||v||_{G} = \sqrt{h(X)}.||v||_{W} \mbox{ for } v \in T_{X}(QF(S)).$$
 \end{theorem}
 
 {\bf Proof:}
From equation  \ref{metricdefalt}, we have
$$||v||^{2}_{W} = \frac{\int \ddot{\Phi}_{0} \ dm}{\int \Phi_{0}\ dm} = \frac{\frac{d^{2}}{dt^{2}}\left(h(X_{t}).\int \phi_{t} \ dm\right)|_{t=0}}{h(X).\int \phi_{0}\ dm} = \frac{1}{h(X)}(h.F)''(v)$$
were
$$F(t) = \frac{\int \phi_{t}\ dm}{\int \phi_{0}\ dm}.$$
Therefore the result follows from showing that $F(t) = L_{\mu_{X}}(X(t))$ where $\mu_{X} \in {\cal C}(S)$ is the unit Patterson-Sullivan geodesic current for $X$. 

By the density of discrete geodesic currents (theorem \ref{dense}), there exists a sequence of discrete geodesic currents $\mu_{n}$ such that $\mu_{n} \rightarrow \mu_{X}$.  As $\mu_{X}$ is unit length in $X$ we can normalize so that $\mu_{n}$ are unit length in $X$. Therefore $\mu_{n} = \a_{n}/l_{n}$ where $\a_{n}$ is a geodesic current coming from Dirac measure on the lifts of a primitive geodesic (also labeled $\a_{n}$) and $l_{n}$ is the length of $\a_{n}$ in $X$.

 We choose as our basepoint for $QF(S)$ the fuchsian group  $\C = \C_{r}$ described in section \ref{Bowensection}. 
 By lemma \ref{techlemma}  for each $\a_{n}$ we can choose a lift $\gamma_{n} \in \C$ such that the axis $g_{n}$ of $\gamma_{n}$  abuts $D$ and has  fixed points $a_{n}, b_{n}$ with expanding fixed point   $a_{n}$ being a periodic point for $f_{\C}$ with period $p_{n}$. 
 
As $f^{p_{n}}_{\C}(a_{n}) = a_{n}$, then $f_{\C}^{p_{n}} = \gamma \in \C$ in an open neighborhood of $a_{n}$ where $\gamma$ fixes $a_{n}$. Then $\gamma$ and $\gamma_{n}$ both fix $a_{n}$  and therefore are comensurate with axes being equal. As $\gamma_{n}$ is primitive, it follows that $\gamma = \gamma^{k_{n}}_{n}$  for some non-zero integer $k_{n}$. Letting $X(t) = [g_{t},\Gamma_{t}]$, then  $g_{t}(a_{n})$ is a fixed  point of $f^{p_{n}}_{\Gamma_{t}}$. Also if we let $\gamma_{n,t} = g_{t} \circ \gamma_{n} \circ g_{t}^{-1}$ then  $\gamma_{n,t}$ has fixed point $g_{t}(a_{n})$ and  $f_{\Gamma_{t}}^{p_{n}} = (\gamma_{n,t})^{k_{n}}$ in an open neighborhood of $g_{t}(a_{n})$.
 
 Then for $i = 0,p_{n}-1$ we let $g_{n}(i)$ be the element of the orbit of $g_{n}$  with endpoints $a_{n}(i)$ and $b_{n}(i)$, where $a_{n}(i) = f^{i}_{\C}(a_{n})$. By lemma \ref{goodgeodesics}, as the geodesic  $g_{n}$ abuts $D$, then $g_{n}(i)$  must also abut $D$. Therefore by compactness of the finite union of domains abutting $D$, there is an $\e > 0$ such that
  $|a_{n}(i)-b_{n}(i)| > \e$ for all $n, i$.

 We let $m_{n}$ be the probability  measure on $\Sph^{1}$  obtained by taking $1/p_{n}$ dirac measure on the $a_{n}(i), i = 0,\ldots,p_{n-1}$.
 
We have $\phi_{t}:\Sigma \rightarrow \Real$ is given by $\phi_{t}(x) = -\log|f_{\Gamma_{t}}'(\pi_{t}(x))|$ where $f_{\Gamma_{t}} = g_{t}\circ f_{\C} \circ g^{-1}_{t}$ and $\pi_{t} = g_{t} \circ \pi_{\C}$.
Therefore $\phi_{t}(x) = \overline{\phi}_{t}(\pi_{\C}(x))$ where $\overline{\phi}:\Sph^{1} \rightarrow \Real$ is  the map $\overline{\phi}_{t}(z) = -\log|f_{\Gamma_{t}}'(g_{t}(z))|$ .
  Then
 $$m_{n}( \overline{\phi}_{t}) = \int_{\Sph^{1}} \overline{\phi}_{t}(z)\  dm_{n} = -\frac{1}{p_{n}} \log|(f^{p_{n}}_{\Gamma_{t}})'(g_{t}(a_{n}))|. $$
 As $g_{t}(a_{n})$ is a  fixed  point of $f^{p_{n}}_{\Gamma_{t}}$ with $f_{\Gamma_{t}}^{p_{n}} = \gamma_{t}^{k_{n}}$ in an open neighborhood of $g_{t}(a_{n})$ we have
$$m_{n}(\overline{\phi}_{t}) =\int_{\Sph^{1}} \overline{\phi}_{t}\  dm_{n}  = -\frac{1}{p_{n}}\log|((\gamma_{n,t})^{k_{n}})'(g_{t}(a_{n}))| = -\frac{k_{n}}{p_{n}} L_{\gamma_{n}}(X_{t}).$$

In particular we have 
$$m_{n}(\overline{\phi}_{0})   = -\frac{k_{n}}{p_{n}} L_{\gamma_{n}}(X) = -\frac{k_{n}l_{n}}{p_{n}}.$$
Therefore
$$\frac{m_{n}(\overline{\phi}_{t})}{m_{n}(\overline{\phi}_{0})} = \frac{L_{\gamma_{n}}(X_{t})}{L_{\gamma_{n}}(X)} = L_{\mu_{n}}(X_{t}).$$
We now show that $l_{n}/p_{n}$ is bounded. 
As the map $\overline{\phi}_{0}$ is bounded on $\Sph^{1}$, there exists a $C$ such that $|\overline{\phi}_{0}| \leq C$. As  $\mu_{n}$ is a probability measure and $k_{n}$ is a non-zero integer,
\begin{equation}
\frac{l_{n}}{p_{n}} \leq \left| \frac{k_{n}l_{n}}{p_{n}} \right| \leq \left|\int \overline{\phi}_{0}\  dm_{n} \right|  \leq \int |\overline{\phi}_{0}| dm_{n} \leq C.
\label{C}
\end{equation}

  Let $\nu_{n}$ be the probability measure on $G(\Hyp^{3})$ obtained by taking $1/p_{n}$ times Dirac measure on set  of geodesics $g_{n}(i)$ given by the endpoint pair ${(a_{n}(i),b_{n}(i))}$. 
As $|a_{n}(i) - b_{n}(i)| > \e$ for all $n,i$,  the measures $\nu_{n}$ does not accumulate on the diagonal,  therefore the sequence $\nu_{n}$ has convergent subsequences in the weak$^{*}$ topology on $G(\Hyp^{3})$. Let $\nu$ be a limit with $\nu = \lim_{i \rightarrow \infty }\nu_{n_{i}}$.  

We will show that $\nu$ is absolutely continuous with respect to  $\mu_{X}$. Let $\mu_{X}(A) = 0$,  then as $\mu_{n} \rightarrow \mu_{X}$
$$\lim_{n \rightarrow \infty} \mu_{n}(A) = \mu_{X}(A) = 0.$$
We compare $\mu_{n}$ and $\nu_{n}$. Both are discrete measures and the support of $\nu_{n}$ is contained in the support of $\mu_{n}$ and with measures $\nu_{n}, \mu_{n}$ having point masses $1/p_{n}, 1/l_{n}$ respectively. Therefore by equation \ref{C}
$$\nu_{n}(A) \leq \frac{l_{n}}{p_{n}} \mu_{n}(A) \leq C. \mu_{n}(A).$$
Thus  we have
$$\nu(A) = \lim_{i \rightarrow \infty} \nu_{n_{i}}(A) \leq \lim_{i \rightarrow \infty} C.\mu_{n_{i}}(A)  \leq C.\mu_{X}(A) = 0.$$
Thus $\mu_{X}(A) = 0$ implies $\nu(A) = 0$. Thus $\nu$ is absolutely continuous with respect to $\mu_{X}$. We take $m_{n_{i}}$ to be the probability measures corresponding to the convergent sequence $\nu_{n_{i}}$. By reducing to subsequence we can assume that $\mu_{n_{i}}$ converge to a  probability measure $m_{f}$ on $\Sph^{1}$. Then we have that $m_{f}$ satisfies
$$m_{f} = \lim_{i  \rightarrow \infty} m_{n_{i}}.$$
and 
$$\frac{m_{f}(\overline{\phi}_{t})}{m_{f}(\overline{\phi}_{0})} = \lim_{i \rightarrow \infty} \frac{m_{n_{i}}(\overline{\phi}_{t})}{m_{n_{i}}(\overline{\phi}_{0})} = \lim_{i \rightarrow \infty} L_{\mu_{n_{i}}}(X_{t}) = L_{\mu_{X}}(X_{t}).$$

Let $g_{0}:\Sph^{2}\rightarrow \Sph^{2}$ be the quasiconfromal homeomorphism conjugating $\C$ to $\Gamma_{0}$. Then by the definition of the Patterson-Sullian geodesic current $\mu_{X}$ (see equation  \ref{psdiff}) we have
$$(d\mu_{X})(a,b) = \pi_{*}\left(\frac{dm_{X}dm_{X}}{|g_{0}(a)-g_{0}(b)|^{2h(X)}}\right).$$
where $m_{X}$ is the Patterson-Sullivan measure for $\Gamma_{0}$ and $\pi$ is the  $\Z_{2}$ cover $\pi: (\Sph^{2}\times \Sph^{2} -\mbox{diagonal}) \rightarrow G(\Hyp^{3})$.
Therefore $\mu_{X}$ is absolutely continuous with respect  to the measure $\pi_{*}(m_{X}\times m_{X})$ on $G(\Hyp^{3})$.

Now we will show that $m_{f}$ is absolutely continuous with respect to $m_{X}$. Let $m_{X}(A) = 0$.
If  $S \subset \Sph^{2}\times \Sph^{2}$, we let $[S] = \pi(S-\mbox{diagonal})$. Then $[S]$ is precisely the set of (unoriented) geodesics  in $S$. Then by definition of $m_{n}$ we have that
$$m_{n}(A) = \nu_{n}([A\times \Sph^{2}]).$$
As $m_{X}(A) = 0$ then $(m_{X} \times m_{X})(A \times \Sph^{2}) = (m_{X} \times m_{X}) (\Sph^{2}  \times A) = 0$. Therefore on $G(\Hyp^{3})$ we obtain
$$\pi_{*} (m_{X} \times m_{X}) ([A\times\Sph^{2}]) = m_{X}\times m_{X}(\pi^{-1}( [A\times \Sph^{2}])) $$
$$= m_{X}\times m_{X}(\left((A\times \Sph^{2})\cup (\Sph^{2}\times A)\right)-\mbox{diagonal})$$
$$  \leq (m_{X} \times m_{X}) (A \times \Sph^{2}) + (m_{X} \times m_{X}) (\Sph^{2}  \times A)  = 0$$
Therefore $\pi_{*} (m_{X} \times m_{X}) ([A\times\Sph^{2}])  = 0$ and as $\mu_{X}$ is absolutely continuous with respect  to $\pi_{*}(m_{X}\times m_{X})$, then $\mu_{X}([A\times \Sph^{2}]) = 0$. Then as $\nu$ is absolutely continuous with respect to $\mu_{X}$ we have $\nu([A \times \Sph^{2}]) = 0.$
As  $m_{n}(A) = \nu_{n} ([A \times \Sph^{2}])$ then
$$m_{f}(A) = \lim_{i \rightarrow \infty} m_{n_{i}}(A) =  \lim_{i \rightarrow \infty} \nu_{n_{i}} ([A \times \Sph^{2}]) = \nu([A \times \Sph^{2}]) = 0.$$
Therefore the limit $m_{f}$ must be absolutely continuous with respect to $m_{X}$.  By Sullivan \cite{Sul79}, the Patterson-Sullivan measure $m_{X}$ is equal to the Hausdorff measure of dimension $h(X)$ on the limit set. Also by Bowen  the pushforward $\overline{m} = (\pi_{\C_{0}})_{*}(m)$ of  the equilibrium measure $m$ to $\Sph^{1}$ is equivalent to the Hausdorff measure of dimension $h(X)$ on the limit set and therefore equivalent to $m_{X}$ ((lemma 10 of \cite{Bow79}). Therefore  $m_{f}$ is then absolutely continuous with respect to the measure $\overline{m}$. Also as the $m_{n}$ are invariant under $f_{\C_{0}}$, then the limit $m_{f}$ is invariant under $f_{\C_{0}}$.  As $\overline{m}$ is ergodic, then $m_{f}$ is also ergodic. But by the Ruelle-Perron-Frobinius theorem (see theorem \ref{rpf}),   there is a unique $f_{\C_{0}}$ invariant ergodic probability measure. Thus $m_{f} = \overline{m}$. 
 As $\phi_{t} = \overline{\phi}_{t} \circ \pi_{\C}$,  then
$m(\phi_{t}) = \overline{m}(\overline{\phi}_{t})$ and
$$F(t) = \frac{m(\phi_{t})}{m(\phi_{0})} = \frac{\overline{m}(\overline{\phi}_{t})}{\overline{m}(\overline{\phi}_{0})} =  \frac{m_{f}(\overline{\phi}_{t})}{m_{f}(\overline{\phi}_{0})} = L_{\mu_{X}}(X_{t}).$$
\eproof

\section{Positive-definite locus for $G$}
Before we prove the main theorem we characterize the zero vectors of $G$ in terms of derivatives of length functions.

\begin{theorem}
Let $v \in T_{X}(QF(S))$ then $||v||_{G} = 0$ if and only if 
for every $\gamma \in \C$, the associated length function $L_{\gamma}:QF(S) \rightarrow \Real$ satisfies
$$(h.L_{\gamma})'(v) = 0.$$
\label{Wandlength}
\end{theorem}

{\bf Proof:}  We  choose our basegroup $\C$ to be the fuchsian group  described in lemma \ref{techlemma}, i.e. if $S$ is a genus $g$ surface, then $\C$ is generated by the standard identification of the sides of the regular $4g$-gon.

We first prove that $||v||_{G} = 0$ implies that $(h.L_{\gamma})'(v) = 0$ for all $\gamma \in \C$.
As it is automatically true for $v = 0$, we assume that $v \neq 0$ and  choose a smooth curve $\a:(-\e,\e) \rightarrow QF(S)$ with $\a(0) = X$ and $\a'(0) = v$ and $\a(t)= X_{t} = [g_{t},\Gamma_{t}]$ as before. Therefore $g_{t}$ conjugates the action of $\C$ to the action of $\Gamma_{t}$.

We now use a trick to reduce the problem to showing that $(h.L_{\gamma})'(v) = 0$  for $\gamma$ a certain subset of $\C$ and then using geodesic currents to show that it is true for all of $\gamma \in \C$. 

We define the map $F:{\cal C}(S) \rightarrow \Real$ by
$$F(\mu) = \frac{(h.L_{\mu})'(v)}{L_{\mu}(X)}.$$ 
Then $F$ is continuous on ${\cal C}(S)$ and constant on positive rays $\{r.\mu\ | r \in \Real_{+}\}$. We define  the space of projective currents 
${\cal PC}(S) = {\cal C}(S)/\sim$ where $\mu_{1} \sim \mu_{2} \mbox{ if }\mu_{2} = r.\mu_{1}$ for some $r \in \Real_{+}$. Then we have the continuous map $\overline{F}:{\cal PC}(S) \rightarrow \Real$ by $\overline{F}([\mu]) = F(\mu)$. Therefore by continuity,  the theorem follows if we prove  $\overline{F} = 0$ on  a dense subset of ${\cal PC}(S)$.
As the set of discrete geodesic currents   are dense in ${\cal C}(S)$ (see theorem \ref{dense}), the set of projective discrete geodesic currents, labeled ${\cal DPC}(S)$, is dense in ${\cal PC}(S)$. Also any set containing all but a finite set of projective discrete geodesic currents is dense in ${\cal PC}(S)$. 

If $\mu \in {\cal C}(S)$ is a discrete geodesic current, then $\mu = r.\a$ where $\a$ is the geodesic current of primitive closed geodesic. We then let $\gamma \in \C$ be an element corresponding to a lift of $\a$. Then we see that the set of projective discrete geodesic currents ${\cal DPC}(S)$ is naturally equivalent to the set of commensurability classes $[\C]$ by the map  $[\mu] \rightarrow [\gamma]$.
 Also we note that by definition the length functions satisfy $L_{\mu} = k.L_{\gamma}$. Thus we need only prove that
$(h.L_{\gamma})'(v) = 0$ for all but a finite set of commensurability classes  $[\gamma]$. We will choose this set to be ${\cal D} = \{[\gamma]\ |,[\gamma] \not\in S\}$, where $S$ is the finite set defined in the above lemma \ref{techlemma}.

We note that if $M$ is a loxodromic m\"obius transformation then the translation distance of $M$ is given by $\log|M'(z)|$ where $z$ is the expanding fixed point of $M$.

 Thus if $\gamma \in  \C$, with expanding fixed point  $z$ then we let 
 $\gamma_{t} = g_{t}\circ \gamma \circ g_{t}^{-1} \in \Gamma_{t}$. Then $
 \gamma_{t}$ has expanding fixed point $z_{t} = g_{t}(z)$ and
 \begin{equation}
 L_{\gamma}(X_{t}) = \log|\gamma'_{t}(z_{t})|.
 \label{lengthequation}
 \end{equation}
 
As $||v||_{G} = 0$, then as $G$ is conformally equivalent  to $W$, $||v||_{W} = 0$.  Therefore by equation \ref{metricdefalt},  $Var(\dot{\Phi}_{0}, m) = 0$. But by non-degeneracy of the variance, this gives
$\dot{\Phi}_{0} \sim 0$ and is a coboundary. Therefore there is a  continuous function $u:\Sigma \rightarrow \Real$ such that
$\dot{\Phi}_{0}(x) = u(\sigma(x)) - u(x).$ Iterating we have 
$$(S_{n}\dot{\Phi}_{0})(x) = u(\sigma^{n}(x)) - u(x).$$
In particular if $\sigma^{n}(x) = x$ then $(S_{n}\dot{\Phi}_{0})(x) = 0$.

Now let $[\gamma] \in {\cal D}$. Then by lemma \ref{techlemma},  there is an element $\gamma'  \in [\gamma]$ with expanding fixed point $z$ is good and a periodic point of $f_{\C}$. Therefore there is an $n$ such that $f_{\C}^{n}(z) = z$. As $z$ is a good point, we let $x \in \Sigma$ be the unique point such that $\pi_{\C}(x) = z$. As $\pi_{\C}\circ \sigma = f_{\C}\circ \pi_{\C}$, we have
$$\pi_{\C}(\sigma^{n}(x)) = f^{n}_{\C}(\pi_{\C}(x)) = f^{n}_{\C}(z) = z = \pi_{\C}(x).$$
As $z$ is a good point, $x$ is the unique preimage of $z$ under $\pi_{\C}$. Therefore  $\sigma^{n}(x) = x$.

As $f_{\C}$ is a Markov map and $z$ is a good point, there is a $\gamma_{z}  \in \C$ such that $f^{n}_{\C} = \gamma_{z}$ in an open neighborhood of $z$. Also as $f_{\C}$ is an expanding Markov map, $z$ is the expanding fixed point of $\gamma_{z}$. Thus elements $\gamma'$ and $\gamma_{z}$ have  common fixed point $z$. As $\C$ is co-compact, this implies that $\gamma'$ and $\gamma_{z}$ are commensurate. Therefore by transitivity of commensurability, $\gamma$ and $\gamma_{z}$ are commensurate.

We let $\gamma_{z,t} = g_{t} \circ \gamma_{z} \circ g_{t}^{-1} \in \Gamma_{t}$. By definition $\pi_{t}(x) = g_{t}(\pi_{\C}(x)) = g_{t}(z)$ and is therefore the expanding fixed point of 
$\gamma_{z,t}$. Also as $f_{\Gamma_{t}} = g_{t} \circ f_{\C}\circ g_{t}^{-1}$ we have
$f^{n}_{\Gamma_{t}} = g_{t} \circ f^{n}_{\C}\circ g_{t}^{-1}$ and therefore
$f^{n}_{\Gamma_{t}} = \gamma_{z,t}$ at $g_{t}(z)$. Thus by the iteration relation in equation \ref{iteration} we have
$$(S_{n}\phi_{t})(x) = -\log|(f^{n}_{\Gamma_{t}})'(\pi_{t}(x))| = -\log|\gamma'_{z,t}(z_{t})| = -L_{\gamma_{z}}(X_{t}).$$
Also  as $\Phi_{t}(x) = h(X_{t})\phi_{t}(x)$, $(S_{n}\Phi_{t})(x) = h(X_{t}).(S_{n}\phi_{t})(x)$. Therefore
$$(S_{n}\dot{\Phi}_{0})(x) = \frac{d}{dt} \left((S_{n}\Phi_{t})(x)\right)\left|_{t=0}  \right.= \frac{d}{dt}\left(-h(X_{t}).L_{\gamma_{z}}(X_{t})\right)\left|_{t=0} = -(h.L_{\gamma_{z}})'(v)\right.$$
As $(S_{n}\dot{\Phi}_{0})(x)  = 0$,   we have that $(h.L_{\gamma_{z}})'(v) = 0$. 
Therefore we have the continuous function $\overline{F}:{\cal PC}(S)  \rightarrow \Real$ is zero on a dense set of points and is therefore the zero function. Thus $(hL_{\mu})'(v) = 0$ for all $\mu \in {\cal C}(S)$ and in particular $(h.L_{\gamma})'(v) = 0$ for all $\gamma \in \C$.

We now prove that if $v$ satisfies $(h.L_{\gamma})'(v) = 0$ for all $ \gamma \in \C$ then $||v||_{G} = 0$.
As it is true for $v = 0$, we assume  that $v \neq 0$ and as before, choose a smooth curve $\a:(-\e,\e) \rightarrow QF(S)$ with $\a(0) = X$ and $\a'(0) = v$ and $\a(t)= X_{t} = [g_{t},\Gamma_{t}]$.

A theorem of Livsic  that states $f \sim g$ if and only if $(S_{n}f)(x) = (S_{n}g)(x)$ whenever $\sigma^{n}(x) = x$ (see \cite{Liv72}). Therefore we let $\sigma^{n}(x) = x$. Then for $z = \pi_{\C}(x)$ we have
$f^{n}_{\C}(z) = z$ and  $f^{n}_{\C} = \gamma_{z}$ at $z$ for some $\gamma_{z} \in \C$. As above we have
$(S_{n}\phi_{t})(x)  = -L_{\gamma_{z}}(X_{t})$ and
$$(S_{n}\dot{\Phi}_{0})(x) = -(h.L_{\gamma_{z}})'(v).$$
 By the assumption $(h.L_{\gamma_{z}})'(v) = 0$ and therefore $(S_{n}\dot{\Phi}_{0})(x) = 0$.
 Therefore by the result of Livsic, $\dot{\Phi}_{0} \sim 0$ and therefore $Var(\dot{\Phi}_{0}, m) = 0$. It follows that $||v||_{W} = 0.$ As $G$ is conformally equivalent to $W$, then $||v||_{G}=0$.
\eproof

 \begin{corollary}
If $||v||_{W} = 0$ then there is a $k \in \Real$ such that
$$L'_{\mu}(v) = k.L_{\mu}(X) \mbox{ for all } \mu \in C(S).$$
 \label{maincorollary}
 \end{corollary}
 {\bf Proof:}
 If $||v||_{W} = 0$  then $(h.L_{\mu})'(v) = 0$ for all $\mu \in C(S)$. Therefore
 $$h'(v)L_{\mu}(X) + h(X).L'_{\mu}(v) = 0.$$
 Solving we have
 $$L'_{\mu}(v) =\left(\frac{-h'(v)}{h(X)}\right).L_{\mu}(X) = k.L_{\mu}(X).$$
\eproof

\medskip

We let $S$ be a closed hyperbolic surface with $S = \Hyp^{2}/\Gamma$ as before and let  $g \in \Gamma$. Given any $X = [f_{0},\Gamma_{0}] \in QF(S)$, then $g$ can be identified to a unique element $g(\Gamma_{0}) = f_{0}\circ g\circ f_{0}^{-1} \in \Gamma_{0} \subseteq PSL(2,\Cplex)$. We can conjugate such that 
$g(\Gamma_{0})$ is of the form 
$$\pm \left( \begin{array}{cc}
\lambda_{g}(X) & 0\\
0 &   \lambda_{g}^{-1}(X)\end{array}\right) \in PSL(2,\Cplex), \mbox{ where } |\lambda_{g}(X)| > 1.$$  We note that $\lambda_{g}$ is well-defined up to sign and $\lambda^{2}_{g}(X)$ is therefore well-defined.

Therefore the element $g(\Gamma_{0})$ is conjugate to  the fractional linear map $f(z) =c.z$, where $c =  \lambda_{g}^{2}(X)$. Therefore we have that the length function $L_{g}:QF(S) \rightarrow \Real$ is given by $L_{g}(X) = 2.\log |\lambda_{g}(X)|$.
Also the holomorphic length function ${\cal L}_{g}:QF(S) \rightarrow \Cplex$ 
satisfies $L_{g} = \Re({\cal L}_{g})$ and $\lambda_{g}^{2} = e^{{\cal L}_{g}}$.

Let $X:(-\e,\e) \rightarrow QF(S)$ be a smooth curve such that $X'(0) = v$. We let $X(t) = [f_{t},\Gamma_{t}]$. Let $\Gamma_{t}$ be a smooth parameterization. Thus  for $g \in \Gamma_{0}$, the map
 $\gamma_{g}:(-\e, \e)  \rightarrow PSL(2,\Cplex)$ defined by $\gamma_{g}(t) = g(\Gamma_{t})$ is  a smooth function. Also as $g(\Gamma_{t}) \in PSL(2,\Cplex)  = SL(2,\Cplex)/\pm I$, we can lift $\gamma_{g}$ to a smooth map $\tilde{\gamma}_{g}:(-\e,\e)  \rightarrow SL(2,\Cplex)$.

We then can define $\lambda_{g}:(-\e,\e) \rightarrow \Cplex$ by letting $\lambda_{g}(t)$ equal the largest eigenvalue of $\tilde{\gamma}_{g}(t)$. Furthermore we define the trace functions
$$t_{g}(t) = tr(\tilde{\gamma}_{g}(t))  = \lambda_{g}(t) + \lambda_{g}^{-1}(t).$$

\begin{lemma}
Let $v \in T_{X}(QF(S))$, $v \neq 0$. If there exists a $k \in \Real$ such that 
$$L'_{g}(v) = k.L_{g}(X) \mbox{ for all } g \in \Gamma$$
then $\lambda^{2}_{g}$,  and $t_{g}^{2}$ are both real  and
$$\Re\left(\frac{\lambda'_{g}}{\lambda_{g}}\right) = 0$$
for all $g \in \Gamma$.
\label{mainlemma}
\end{lemma}

{\bf Proof:} As trace functions are holomorphic co-ordinate function for $QF(S)$ (see \cite{Mar74}), as $v \neq 0$, there exists  $\a \in \Gamma$ be such that $t'_{\a}(0) \neq 0$.
As
$$t_{g}' = \lambda_{g}'. - \frac{1}{\lambda_{g}^{2}}. \lambda_{g}' = \lambda_{g}'.\left(\frac{\lambda_{g}^{2}-1}{\lambda_{g}^{2}}\right)$$
 then $\lambda_{\a}'(0) \neq 0$.

As $\Gamma$ is non-elementary, we can choose a $\b \in \Gamma$ such that $\a,\b$ do not have the same axis. We note that $\a,\b$ have the same axes if and only if there exist $n,m \in \Z$, both non-zero, such that $\a^{n}= \b^{m}$.

By conjugation of $\Gamma_{t}$ we can put $\a(\Gamma_{t})$ in the diagonal form with 

$$A(t) = \tilde{\gamma}_{\a}(t)  = \left( \begin{array}{cc}
\lambda_{\a}(t) & 0\\
0 &   \lambda_{\a}^{-1}(t)\end{array}\right) 
$$ where $|\lambda_{\a}(t)| > 1$.
Therefore we have that
 
$$B(t) = \tilde{\gamma}_{\b}(t)  = \left( \begin{array}{cc}
a(t) & b(t)\\
c(t) &  d(t) \end{array}\right).$$
where $a(t).d(t)-b(t).c(t) = 1.$

We consider the two generator subgroup $G_{t} = <A(t),B(t)> \subseteq SL(2, \Cplex)$  acting on upper half space by the associated fractional linear maps. Then $A(t)$ fixes $0, \infty$ and has axis the $z-$axis. If $a(t) = 0$ then $B(t)$ sends $\infty$ to $0$  and if  $d(t) = 0$ then $B(t)$ sends $0$ to $\infty$. In either case $C(t) = B(t)A(t)B(t)^{-1}$ fixes either $0$ or $\infty$.  As $C(t)$ and $A(t)$ share a fixed point, and  $\Gamma_{t}$ has no parabolics, $C(t)$ and $A(t)$ must have the same axes. Thus $B(t)$ must send both the point $0$ to $\infty$ and $\infty$ to $0$ and has the same axis as $A(t)$ giving our contradiction. Therefore we have that $a(t), d(t)$ are  both non-zero.
  
If $||v||_{W} = 0$, then by corollary \ref{maincorollary}, $L_{g}'(v) = k.L_{g}(X)$ for all $g \in \Gamma$. As $L_{g} = \log|\lambda_{g}|$, we obtain the equation
\begin{equation}
(\log|\lambda_{g}|)'(v) = k.\log|\lambda_{g}(X)|.
\end{equation} 

As we are only interested in derivatives at $0$ for $X(t)$, we will make the notation that $f' = f'(0)$.

Therefore for $g = \alpha$ we have $(\log|\lambda_{\a}|)' = k.\log|\lambda_{\a}|$
or equivalently
\begin{equation}
(\log|\lambda_{\a}|)' =  \frac{|\lambda_{\a}|'}{|\lambda_{\a}|}  = k.\log|\lambda_{\a}|.
\label{loglambda}
\end{equation}

Now we consider the element $C_{n} = A^{n}.B$. Then
$$C_{n} = \left( \begin{array}{cc}
\lambda_{\a}^{n}a & \lambda_{\a}^{n}b\\
\lambda_{\a}^{-n}c &  \lambda_{\a}^{-n}d \end{array}\right).$$
Let $\mu_{n}, \mu_{n}^{-1}$ be the eigenvalues of $C_{n}$, with $|\mu_{n}| > 1$ and define $t_{n} = Trace(C_{n}) = \lambda_{\a}^{n}a +  \lambda_{\a}^{-n}d$.  Then we have
$$t_{n} = \lambda_{\a}^{n}a +  \lambda_{\a}^{-n}d = \mu_{n} + \mu_{n}^{-1}.$$

Solving for $\mu_{n}$ we have
$$\mu_{n} = \frac{t_{n} \pm \sqrt{t_{n}^{2}-4}}{2} $$

Expanding out we get
 $$\sqrt{t_{n}^{2}-4}= \sqrt{\lambda_{\a}^{2n}a^{2} +2ad+  \lambda_{\a}^{-2n}d^{2}-4} = \lambda^{n}a. \sqrt{1 + \lambda_{\a}^{-2n}.\left(\frac{2ad-4}{a^{2}}\right) +  \lambda_{\a}^{-4n}\frac{d^{2}}{a^{2}}}.$$
Therefore, for $n$ large positive, we have
$$\sqrt{t_{n}^{2}-4} = \lambda^{n}a. \left(1 + \lambda_{\a}^{-2n}.\left(\frac{ad-2}{a^{2}}\right) + O(\lambda_{\a}^{-4n})\right)$$
and
$$\mu_{n} =  \frac{ \lambda_{\a}^{n}a +  \lambda_{\a}^{-n}d + \sqrt{( \lambda_{\a}^{n}a +  \lambda_{\a}^{-n}d)^{2}-4}}{2} $$  
$$ =  \frac{\lambda_{\a}^{n}a \left(1 +  \lambda_{\a}^{-2n}\frac{d}{a}\right) + \lambda^{n}a. \left(1 + \lambda_{\a}^{-2n}.\left(\frac{ad-2}{a^{2}}\right) + O(\lambda_{\a}^{-4n})\right)}{2}$$
Giving 
\begin{equation}
\mu_{n} =  \lambda_{\a}^{n}a \left(1 +  \lambda_{\a}^{-2n}\left(\frac{ad-1}{a^{2}}\right) + O(\lambda_{\a}^{-4n})\right)
\label{eexpansion}
\end{equation}
Thus for element $g = \alpha^{n}\b$ we have 
$$|\lambda_{g}| = |\mu_{n}| = |\lambda_{\a}|^{n}.|a|\left|  1 +  \lambda_{\a}^{-2n}\left(\frac{ad-1}{a^{2}}\right) + O(\lambda_{\a}^{-4n})\right|.$$
Taking logs we have
$$\log|\lambda_{g}| = n\log|\lambda_{\a}| + \log|a| + \log\left|  1 +  \lambda_{\a}^{-2n}\left(\frac{ad-1}{a^{2}}\right) + O(\lambda_{\a}^{-4n})\right|.$$
Expanding we have
$$\log|\lambda_{g}| = n\log|\lambda_{\a}| + \log|a| +   \Re\left( \lambda_{\a}^{-2n}\left(\frac{ad-1}{a^{2}}\right)\right) + O(|\lambda_{\a}|^{-4n}).$$
Then differentiating
$$(\log|\lambda_{g}|)' = n\frac{|\lambda_{\a}|'}{|\lambda_{\a}|} + \frac{|a|'}{|a|} +   \Re\left(-2n\lambda_{\a}^{-2n-1}\lambda_{\a}'\left(\frac{ad-1}{a^{2}}\right)\right) + \Re\left(\lambda_{\a}^{-2n}\left(\frac{ad-1}{a^{2}}\right)'\right)+ O(|\lambda_{\a}|^{-4n}).$$
By assumption $(\log|\lambda_{g}|)' - k.\log|\lambda_{g}| = 0$.
Therefore for large positive $n$
$$0 = n.\left(\frac{|\lambda_{\a}|'}{|\lambda_{\a}|}-k\log|\lambda_{\a}|\right) + \left(\frac{|a|'}{|a|}-k \log|a|\right) +   \Re\left(-2n\lambda_{\a}^{-2n-1}\lambda_{\a}'\left(\frac{ad-1}{a^{2}}\right)\right) +$$
$$ \Re\left(\lambda_{\a}^{-2n}\left(\left(\frac{ad-1}{a^{2}}\right)'  -k.\left(\frac{ad-1}{a^{2}}\right)\right)\right) +  O(|\lambda_{\a}|^{-4n}).$$
We now derive the equations we are looking for. Taking limits we have
$$\lim_{n \rightarrow \infty}\frac{(\log|\lambda_{g}|)' - k.\log|\lambda_{g}|}{n} = \frac{|\lambda_{\a}|'}{|\lambda_{\a}|} - k.\log|\lambda_{\a}| = 0.$$
This is just equation \ref{loglambda} we already obtained. Taking further limits we have
$$\lim_{n \rightarrow \infty}\left((\log|\lambda_{g}|)' - k.\log|\lambda_{g}|\right) = \frac{|a|'}{|a|} - k.\log|a| = 0.$$
This gives us a new equation
\begin{equation}
 \frac{|a|'}{|a|}  = k.\log|a|.
 \label{aequation}
 \end{equation}
 Now we take the following limit
$$\lim_{n \rightarrow \infty}\frac{|\lambda_{\a}|^{2n}\left((\log|\lambda_{g}|)' - k.\log|\lambda_{g}|\right)}{n}
=  \lim_{n\rightarrow \infty} \Re\left(-2|\lambda|^{2n}|.\lambda_{\a}^{-2n-1}\lambda_{\a}'\left(\frac{ad-1}{a^{2}}\right)\right) = 0$$
Simplifying we get
\begin{equation}
 \lim_{n\rightarrow \infty} \Re\left(\left(\frac{\lambda_{\a}}{|\lambda_{\a}|}\right)^{-2n}.\left(\frac{\lambda_{\a}'}{\lambda_{\a}}\right)\left(\frac{ad-1}{a^{2}}\right)\right) = 0
 \end{equation}
 We let 
$$u = \left(\frac{\lambda_{\a}}{|\lambda_{\a}|}\right)^{2}.$$
As we can always choose a sequence $n_{i}$ such that $\lim_{i \rightarrow \infty} u^{-n_{i}} = 1$, we have that
 $$\lim_{i \rightarrow \infty} \Re\left(\left(\frac{\lambda_{\a}}{|\lambda_{\a}|}\right)^{-2n_{i}}.\left(\frac{\lambda_{\a}'}{\lambda_{\a}}\right)\left(\frac{ad-1}{a^{2}}\right)\right) = \Re\left(\left(\frac{\lambda_{\a}'}{\lambda_{\a}}\right)\left(\frac{ad-1}{a^{2}}\right)\right) = 0.$$
 Therefore we obtain the equation
\begin{equation} \Re\left(\left(\frac{\lambda_{\a}'}{\lambda_{\a}}\right)\left(\frac{ad-1}{a^{2}}\right)\right) = 0.
\label{realequation}
 \end{equation}
If $\lambda_{\a}^{2}$ is not real, then we let $u = e^{\pi i  \theta}$ for $\theta \in [0,2)$.

{\bf Case 1: $\theta$ irrational:} 
If $\theta$ is irrational, then we can choose a sequence $m_{i}$ such that $\lim_{i \rightarrow \infty} u^{-m_{i}} = i.$
Then
$$\lim_{i \rightarrow \infty} \Re\left(\left(\frac{\lambda_{\a}}{|\lambda_{\a}|}\right)^{-2m_{i}}.\left(\frac{\lambda_{\a}'}{\lambda_{\a}}\right)\left(\frac{ad-1}{a^{2}}\right)\right) = \Im\left(\left(\frac{\lambda_{\a}'}{\lambda_{\a}}\right)\left(\frac{ad-1}{a^{2}}\right)\right) = 0 $$
Thus both the real and imaginary parts are zero giving
$$\left(\frac{\lambda_{\a}'}{\lambda_{\a}}\right)\left(\frac{ad-1}{a^{2}}\right) = 0$$
As $\lambda_{\a}' \neq 0$ we have $ad=1$. Therefore as $ad-bc = 1$, we have $bc = 0$ and either $b = 0$ or $c = 0$. If $b = 0$, then $\a,\b$ have common fixed point $0$ and if $c=0$, then $\a,\b$ have common fixed point $\infty$. As they do not have common fixed points, we have that $u$ is not irrational.

{\bf Case 2: $\theta$ positive rational but not integer.} We let $\theta = p/q$, where $q > 1$ and $p,q$ have no common divisors. Then $u^{q} = 1$ and $u^{nq + 1} = u$.
Then let $n_{i} = i.q - 1$. Then $u^{-n_{i}} = u$. Thus
$$\lim_{i \rightarrow \infty} \Re\left(\left(\frac{\lambda_{\a}}{|\lambda_{\a}|}\right)^{-2n_{i}}.\left(\frac{\lambda_{\a}'}{\lambda_{\a}}\right)\left(\frac{ad-1}{a^{2}}\right)\right) =  \Re\left( u.\left(\frac{\lambda_{\a}'}{\lambda_{\a}}\right)\left(\frac{ad-1}{a^{2}}\right)\right) = 0$$
Let $u = x + iy$ where $y \neq 0$. Then
$$\Re\left( u.\left(\frac{\lambda_{\a}'}{\lambda_{\a}}\right)\left(\frac{ad-1}{a^{2}}\right)\right) = x.\Re\left(\left(\frac{\lambda_{\a}'}{\lambda_{\a}}\right)\left(\frac{ad-1}{a^{2}}\right)\right) - y. \Im\left(\left(\frac{\lambda_{\a}'}{\lambda_{\a}}\right)\left(\frac{ad-1}{a^{2}}\right)\right) = 0$$
Therefore by equation \ref{realequation}, we have
$$y. \Im\left(\left(\frac{\lambda_{\a}'}{\lambda_{\a}}\right)\left(\frac{ad-1}{a^{2}}\right)\right) = 0$$
As $y \neq 0$ we obtain the conclusion that both real and imaginary parts are zero giving
$$\left(\frac{\lambda_{\a}'}{\lambda_{\a}}\right)\left(\frac{ad-1}{a^{2}}\right) = 0$$
This leads to the same contradiction as case 1.

Thus we conclude that $\lambda_{\a}^{2}$ is real. Thus  $\lambda_{\a}$  is either purely imaginary or purely real and as $t_{\a} = \lambda_{\a} + \lambda_{\a}^{-1}$, then $t_{\a}$ is similarly 
either purely imaginary or purely real and $t^{2}_{\a}$ is real. Therefore we have shown that if $t'_{g} \neq 0$ then $t^{2}_{g}$ is real.

Also as $t_{n} = \lambda^{n}_{\a}.a + \lambda^{n}.d$ then
$$t'_{n} = n.\lambda_{\a}^{n-1}\lambda'_{\a}.a + \lambda^{n}_{\a}.a' -  n\lambda_{\a}^{-n-1}\lambda'_{\a}.a + \lambda^{-n}_{\a}.a' .$$
Thus
$$\lim_{n \rightarrow \infty}
\left(\frac{t'_{n}}{n.\lambda^{n}_{\a}}\right) =\frac{\lambda'_{\a}}{\lambda_{\a}}.a$$
and therefore for large $n$, $t'_{n} \neq 0$. Choose $n_{0}$ such that  $t'_{n} \neq 0$ for $n > n_{0}$. 

We let $n > n_{0}$. By the above, $t^{2}_{n}$ is real  and 
$$ t_{n}^{2} = (\lambda_{\a}^{n}a + \lambda_{\a}^{-n}d)^{2} = \lambda_{\a}^{2n}a^{2} + 2ad + \lambda_{\a}^{-2n}d^{2}.$$
 As $t_{n}^{2}$ is real and $\lambda_{\a}^{2}$ is real, we have 
 $$\Im(t^{2}_{n}) =  0 = \lambda_{\a}^{2n}\Im(a^{2}) + 2\Im(ad) + \lambda_{\a}^{-2n}.\Im(d^{2})$$
 Taking limits we have
 $$\lim_{n \rightarrow \infty} \frac{\Im(t_{n}^{2})}{\lambda_{\a}^{2n}} = \Im(a^{2})  = 0.$$
 Therefore
 $$\lim_{n \rightarrow \infty} \Im(t_{n}^{2}) = 2\Im(ad) = 0.$$
  and finally
    $$\lim_{n \rightarrow \infty} (\lambda_{\a}^{2n}.\Im(t_{n}^{2})) = \Im(d^{2}) = 0.$$
    Thus $a^{2}, d^{2}, ad$ are all real. Applying this to  equation \ref{realequation} we have
$$ \Re\left(\left(\frac{\lambda_{\a}'}{\lambda_{\a}}\right)\left(\frac{ad-1}{a^{2}}\right)\right) = \left(\frac{ad-1}{a^{2}}\right).\Re\left(\frac{\lambda_{\a}'}{\lambda_{\a}}\right) = 0.$$
Therefore we have
$$\Re\left(\frac{\lambda_{\a}'}{\lambda_{\a}}\right) = 0.$$
As the only assumption on $\a$ was that $t'_{\a}$ and therefore $\lambda'_{\a}$ was non-zero, we have
$$\Re\left(\frac{\lambda_{g}'}{\lambda_{g}}\right) = 0 \mbox{ for all } g \in \Gamma.$$

    Also as $t^{2}_{\b} = (a + d)^{2} = a^{2} + 2ad + d^{2}$, then we have that $t^{2}_{\b}$ is real. 
    As $\b$ was arbitrarily chosen, we therefore have that $t^{2}_{g}$ is real for all $g \in \Gamma$. As $t^{2}_{g}$ is real, then $\lambda_{g}^{2}$ is also real. 
    \eproof
    
\begin{lemma}
If $v \in T_{X}(QF(S)$ and there exists a $k \in \Real$  such that $L'_{g}(v) = k.L_{g}(X)$ then $k = 0$.
\label{k=0}
\end{lemma}

If $v = 0$ then $L'_{g}(v) = 0$ and obviously $k = 0$. 

Therefore we assume $v \neq 0$.  Let $g \in \Gamma$. We let $\lambda_{g} = |\lambda_{g}|e^{i\theta}$ then
$$\frac{\lambda_{g}'}{\lambda_{g}} = \frac{|\lambda_{g}|'e^{i\theta} + |\lambda_{g}|e^{i\theta}.i\theta'}{|\lambda_{g}|e^{i\theta}} = \frac{|\lambda_{g}|'}{|\lambda_{g}|} + i.\theta'.$$
Thus 
\begin{equation}
\Re\left(\frac{\lambda_{g}'}{\lambda_{g}} \right) =  \frac{|\lambda_{g}|'}{|\lambda_{g}|}
\end{equation}
Then by equation \ref{loglambda}
$$\Re\left(\frac{\lambda_{g}'}{\lambda_{g}} \right) = \frac{|\lambda_{g}|'}{|\lambda_{g}|} = k.\log|\lambda_{g}|.$$

But by the above lemma \ref{mainlemma}
$$\Re\left(\frac{\lambda_{g}'}{\lambda_{g}} \right) = 0.$$
Thus we have $k. \log|\lambda_{\a}| = 0$. As $|\lambda_{\a}| > 1$, $\log|\lambda_{\a}| \neq 0$ and therefore $k = 0$.
\eproof

\begin{lemma}
If $v \in T_{X}(QF(S))$, $v \neq 0$  satisfies
$$L'_{g}(v) = 0,\mbox{ for all } g \in \Gamma$$
 then $X \in F(S)$ and $v = J.w$ for some $w \in T_{X}(F(S)) \subseteq T_{X}(QF(S))$.
\label{lastlemma}
\end{lemma}

{\bf Proof:}
We pick $\a,\b$ as in lemma \ref{mainlemma}. 
For group $G_{0} = <A(0),B(0)> $ we have  $t^{2}_{\a}, \lambda_{\a}^{2}, a^{2}, ad, d^{2}$ are all real. Therefore the fractional linear map given by $A$ is $f_{A}(z) = \lambda_{\a}^{2}.z \in PSL(2,\Real)$. 

As $ad-bc = 1$, we therefore have that $bc = ad-1$ is real. Therefore $b = r.e^{i\theta}$ and $c = s.e^{-i\theta}$ where $r,s$ are real. 

If $a,d$ are both real, we conjugate $G_{0}$ by rotation $R$ about the axis of $A$ by angle $\theta$.
Then as $R,A$ commute, $RAR^{-1}=A$ and
$$RBR^{-1} = 
\left( \begin{array}{cc}
e^{-i\theta/2} & 0\\
0 & e^{i\theta/2}   \end{array}\right)
\left( \begin{array}{cc}
a & b\\
c & d   \end{array}\right)
\left( \begin{array}{cc}
e^{i\theta/2} & 0\\
0 & e^{-i\theta/2}   \end{array}\right)
=
\left( \begin{array}{cc}
a & r\\
s & d   \end{array}\right)
.$$
Therefore the fractional linear map given by $RBR^{-1}$ is  in $PSL(2,\Real)$.

If $a,d$ are both imaginary, we conjugate by a rotation $R$ about the axis of $A$ by angle $\pi+\theta$. Then
$$RBR^{-1} =
\left( \begin{array}{cc}
a & ir\\
is & d   \end{array}\right)
.$$
Thus as each entry is imaginary, the fractional linear map is in $PSL(2,\Real)$.

Therefore we have conjugated $G_{0}$ to a subgroup of $PSL(2,\Real)$. Thus $G_{0}$ has limit set contained in a Euclidean line $L_{0}$ through the origin  and $G_{0}$ preserves a hyperbolic plane $H_{0}$ containing the axis of $A$.  We conclude that if $\a \in \Gamma_{0}$ has $\lambda_{\a}' \neq 0$ then for any $\b \in \Gamma_{0}$, the axes of $\a$ and $\b$ are contained in the same geometric circle.

If $X \not\in F(S)$ then there is an element $\gamma \in \C$ such that the associated fractional linear map $C \in \Gamma_{0}$ does not preserve $H_{0}$.
Then we have as before that the group $G_{1} = <A , C>$ can be conjugated to a subgroup of $PSL(2,\Real)$. Therefore $\Gamma_{1}$ preserves a line $L_{1}$, and hyperbolic plane $H_{1}$ containing the axis of $A$. As by assumption $C$ does not preserve $H_{0}$ then $H_{0} \neq H_{1}$ and therefore $L_{1} \neq L_{0}$. Thus $L_{1} \cap L_{0} = \{0,\infty\}$.
By conjugation, we assume that $L_{0}$ is the real axis. 

We note that if $g,h$ are loxodromic hyperbolic elements, then the axis of $ghg^{-1}$ is the image of the axis of $h$ under $g$.

Thus we conjugate $\a$ by $\b\a^{n}$ to get 
$$\a_{n} = (\b\a^{n})^{-1}\a(\b\a^{n}) = \a^{-n}(\b^{-1}\a\b)\a^{n}$$ 
Then as $\a_{n}$ is a conjugate of $\a$ we have $\lambda'_{\a_{n}} = \lambda_{\a}' \neq 0$.
 We let
$A_{n} = (BA^{n})^{-1})A(BA^{n}) = A^{-n}(B^{-1}AB)A^{n} \in \Gamma_{0}$.  
Then the endpoints of the axes of  $A_{n}$ and $C$ must be contained in a geometric circle.
Also the axis of $A_{n}$ is the image of the axis of $BAB^{-1}$ under $A^{-n}$. Therefore we let $a,b$ be the endpoints of the axis of $B^{-1}AB$. As $A, B$ are non-commensurate, their axes do not have common endpoints. Therefore $a,b \in \Real$ and are not equal $0$ or $\infty$.  Then the endpoint of the axis of $A_{n}$ are  $a_{n}, b_{n}$ where   $a_{n} = A^{-n}(a) =  \lambda_{\a}^{-2n}.a, b_{n} = A^{-n}(b) = \lambda_{\a}^{-2n}.b$. 

Let $z,w \in L_{1}$ be the endpoints of the axis of $C$. As $L_{1}$ is not the real axis, then
$z = re^{i\theta}, w = se^{i\theta}$ where $r,s \in \Real$, are $e^{i\theta}$ is not real.
As the axes of $C$ and $A_{n}$ are on the same geometric circle,  the cross ratio $(a_{n},z;b_{n},w)$ is real for all n.
$$(a_{n},z;b_{n},w) =  \frac{(a_{n}-b_{n})(z-w)}{(a_{n}-w)(z-b_{n})} = 
 \frac{\lambda_{\a}^{-2n}(a-b)(r-s)e^{i\theta}}{(\lambda_{\a}^{-2n}a-se^{i\theta})(re^{i\theta}-\lambda_{\a}^{-2n}b)}$$
Therefore as $\Im (a_{n},z;b_{n},w)  = 0$ and $\lambda_{\a}^{2}$ is real
 then
 $$0 = \lim_{n \rightarrow \infty} \lambda_{\a}^{2n}. \Im (a_{n},z;b_{n},w) = \lim_{n \rightarrow \infty} \Im\left( \frac{(a-b)(r-s)e^{i\theta}}{(\lambda_{\a}^{-2n}a-se^{i\theta})(re^{i\theta}-\lambda_{\a}^{-2n}b)}\right)
 $$
 $$= \Im\left( \frac{(a-b)(r-s)e^{i\theta}}{-sre^{2i\theta}}\right) = \frac{(a-b)(r-s)}{-rs} \Im (e^{-i\theta}) .$$
Thus $ \Im(e^{-i\theta}) = 0$, and therefore $e^{i\theta}$ is real. But by assumption $e^{i\theta}$ which gives us  our contradiction.
Thus $X \in F(S)$.

Finally as $X \in F(S)$, we have the decomposition (see \cite{BT08}),    
$$T_{X}(QF(S)) = T_{X}(F(S)) \oplus J.T_{X}(F(S)).$$ 
 If $v \in T_{X}(F(S))$ then ${\cal L}_{\a}'(v) = L_{\a}'(v)$ and is real. Therefore if $v \in T_{X}(QF(S))$, then $v = v_{1} + J.v_{2}$ where $v_{i} \in T_{X}(F(S))$. Therefore 
$$L'_{g}(v) = \Re({\cal L}_{g}'(v)) = \Re( {\cal L}_{g}'(v_{1}) + {\cal L}'_{g}(J.v_{2})) $$
$$=  \Re({\cal L}'_{g}(v_{1}) + i.{\cal L}'_{g}(v_{2})) = \Re(L'_{g}(v_{1}) + i.L'_{g}(v_{2})) = L'_{g}(v_{1}).$$
Therefore if  $L'_{g}(v) = 0$ for all $g \in \C$, then $L'_{g}(v) = L'_{g}(v_{1}) = 0$ for all $g \in \C$. But this implies that $v_{1} = 0$. Therefore $v = J.v_{2}$
\eproof

We now are ready to prove the main theorem.

\smallskip

\noindent{\bf Proof of main theorem:}

We first prove that if $v \in T_{X}(QF(S)), v \neq 0$, and $||v||_{G} = 0$ then
$X \in F(S)$ and $v = J.w$ for some $w \in T_{X}(F(S))$. 

Let $v \in T_{X}(QF(S)), v \neq 0$, and $||v||_{G} = 0$. As $W$ is a multiple of $G$, $||v||_{W} = 0$. Then by corollary \ref{maincorollary}, there is a $k \in \Real$ such that
$L'_{g}(v) =k.L_{g}(X)$ for all $g \in \Gamma$. Then by lemma \ref{k=0}, $k = 0$. Therefore
 $L'_{g}(v) = 0$ for all $g \in \Gamma$. Finally by lemma \ref{lastlemma}, $X \in F(S)$ and $v = J.w$ for some $w \in T_{X}(F(S))$. 
 
 We now prove that if  $v = J.w$ where $w \in T_{X}(F(S))$ then $||v||_{G} = 0$. By theorem \ref{Wandlength}, we only need to prove that $(hL_{\mu})'(v) = 0$ for all $\mu \in {\cal C}(S)$.
 
As the complex length functions ${\cal L}_{g}$ are holomorphic on $QF(S)$,
 $${\cal L}'_{g}(v) = {\cal L}'_{g}(J.w) = i.{\cal L}'_{g}(w).$$
 As $w \in T_{X}(F(S))$, ${\cal L}'_{g}(w)$ is real and equal ${\cal L}'_{g}(w) = L'_{g}(w)$.
 Therefore ${\cal L}'_{g}(v) = i.L'_{g}(w)$ is purely imaginary giving 
 $$L'_{g}(v) = \Re\left( i.L'_{g}(w) \right) = 0.$$
Thus $L'_{g}(v) = 0$ for all $g \in \Gamma$. As $h$ is minimum on the fuchsian locus $F(S)$ then $h'(v) = 0$ and 
$$(h.L_{g})'(v) = h'(v)L_{g}(X) + h(X).L'_{g}(v) = 0.$$
 Therefore by theorem \ref{Wandlength}, $||v||_{W} = 0.$ As $G$ is conformally equivalent  to $W$ we therefore have $||v||_{G} = 0.$
 \eproof

 \section{Critical points of  Hausdorff dimension}
 We will now use the description  of the positive definite locus of $G$ to obtain information about the critical points of $h:QF(S) \rightarrow \Real$. 
  
 If $f:X \rightarrow \Real$ is a smooth map, then $x \in X$ is a critical point the differential $f'(x): T_{x}(X) \rightarrow \Real$ is the trivial linear function.
 
 If $x$ is a critical point of $f$ then the Hessian of $f$  at $x$ is a well-defined two-form which we label $f''(x)$. If $\dim(X) = n$, the Hessian is a symmetric bilinear form on $T_{x}(X) = \Real^{n}$. Then the {\em signature} of $f''(x)$ is the (well-defined) triple of non-negative integers $(r, s, t), r+s+t =n$, such that there are local co-ordinates $(x_{1}, \ldots, x_{n})$ with 
$$||v||^{2} = \left|\left|\sum_{i=1}^{n} v_{i}\frac{\partial}{\partial x_{i}}\right|\right|^{2} = v^{2}_{1}+ \ldots v^{2}_{r} - v^{2}_{r+1}\ldots - v^{2}_{r+s}.$$ 
We say $f''(x)$ has positive definite of dimension $r$, and negative definite dimension $s$ and trivial dimension $t$.

As $h \geq 1$ and  $h = 1$ on the fuchsian subspace $F(S)$ it follows that each $h$ is minimum (and therefore critical) at each point of  $F(S)$. Thus for $X \in F(S)$, $h''(X)$ has negative definite dimension zero and trivial dimension at least $\dim(F(S)) = 6g-6$. In \cite{BT08}, we show that $h''(X)$ has positive definite dimension $6g-6$. We generalize this to all critical points of $h$ to prove theorem \ref{criticalpoints}.

\noindent{\bf Theorem \ref{criticalpoints}}\newline
{\em If $X \in QF(S)$ is a critical point of $h:QF(S) \rightarrow \Real$ then $X$ has positive definite dimension at least $6g-6$. In particular $h$ has no local maxima.
}

{\bf Proof:}
As the theorem is true for $X \in F(S)$ (see \cite{BT08}), we assume that $X \notin F(S)$. By \cite{BT05}, if $\mu_{X}$ is the unit Patterson-Sullivan geodesic current for $X \in QF(S)$ then  the real valued  function $Y \rightarrow h(Y).L_{\mu_{X}}(Y)$ on $QF(S)$ is minimum at $X$. Therefore $(h.L_{\mu_{X}})'(X) = 0$.

If $X$ is  a critical point of $h$ then $h'(X) = 0$ and therefore by the product rule
$$h'(X)L_{\mu_{X}}(X) + h(X)L'_{\mu_{X}}(X) = h(X)L'_{\mu_{X}}(X) = 0.$$
As $h(X) \neq 0$ then $L'_{\mu_{X}}(X) = 0$ and therefore $L_{\mu_{X}}$ has a critical point at $X$.
We note that the holomorphic length function ${\cal L}_{\mu_{X}}$ satisfies
$$\Re\left({\cal L}_{\mu_{X}} \right) = L_{\mu_{X}}.$$
Therefore as $L'_{\mu_{X}}(X) = 0$, then for all $v \in T_{X}(QF(S))$,
$$\Re\left({\cal L}'_{\mu_{X}}(v) \right) = L'_{\mu_{X}}(v) = 0.$$
Therefore applying this to $J.v$  we have
$$0 = \Re\left({\cal L}'_{\mu_{X}}(J.v) \right) = \Re\left(i.{\cal L}'_{\mu_{X}}(v) \right) =  -\Im\left({\cal L}'_{\mu_{X}}(v) \right).$$
Thus ${\cal L}'_{\mu_{X}}(v)$ has real and imaginary part zero and therefore ${\cal L}'_{\mu_{X}}(v) = 0$ for all $v \in T_{X}(QF(S))$. Thus ${\cal L}'_{\mu_{X}}(X) = 0$ and we have a well-defined complex bilinear 2-form
${\cal L}''_{\mu_{X}}(X).$

As the two-form $G_{X}$ is given by $G_{X} = (hL_{\mu_{X}})''(X)$ we have
$$G_{X}= h''(X)L_{\mu}(X) + 2h'(X)L'_{\mu_{X}}(X)+ h(X)L''_{\mu}(X).$$ 
Therefore as $h$ and  $L_{\mu_{X}}$ are critical at $X$ 
\begin{equation}
G_{X} = h''(X) + h(X)L''_{\mu_{X}}(X).
\label{hsum} 
\end{equation}

Let $L''_{\mu_{X}}$ be positive definite on a subspace $V \subseteq T_{X}(QF(S))$ with $\dim(V) = k$. We now consider the subspace $W = J.V$. If $w \in W$ then $w = J.v$  and 
$$L''_{\mu_{X}}(J.v, J.v) = \Re\left({\cal L}_{\mu_{X}}''(J.v, J.v)\right) =  \Re\left(i^{2}{\cal L}_{\mu_{X}}''(v, v)\right) =  \Re\left(-{\cal L}_{\mu_{X}}''(v, v)\right) = -L''_{\mu_{X}}(v, v).$$
Thus $L''_{\mu_{X}}(X)$ is negative definite on $W$. Therefore $V \cap W = \emptyset$ and $2k \leq \dim(QF(S)) = 12g-12$. Thus $k \leq 6g-6$ and $L''_{\mu_{X}}$ is non-positive on a subspace $V_{1}$ of dimension 
$$\dim(V_{1}) = \dim(QF(S)) - \dim(V)  = (12g-12) - k \geq 6g-6.$$

 As $X \notin F(S)$ then $G_{X}$ is positive definite and 
$$G_{X} = h''(X) + h(X)L''_{\mu_{X}}(X).$$ 
As $L''_{\mu_{X}}(X)$ is non-positive on a subspace $V_{1}$ of dimension at least $6g-6$, then $h''(X)$ must be positive-definite on this subspace. Therefore $h$ has positive definite dimension at least $6g-6$ at $X$.
\eproof
  
  We now give a proof of McMullen's result (theorem  \ref{H=g}) in terms of the description of the positive definite locus of $G$.
  
 {\bf Proof:}
 If $w  = J.v$ for $v \in T_{X}(F(S))$ then by the main theorem we have $||w||_{G} = 0$.
 As $h(X) = 1$, and by holomorphicity, $L''_{\mu_{X}}(J.v, J.v) = - L''_{\mu_{X}}(v)$ we have   
 $$0 = ||w||^{2}_{G} = h''(J.v, J.v) + h(X).L''_{\mu_{X}}(J.v, J.v) = ||v||^{2}_{H}   -L''_{\mu_{X}}(v, v).$$
  Thus
  $$||v||^{2}_{H} = L''_{\mu_{X}}(v, v).$$
   In  \cite{Wol86}, Wolpert  describes the Weil-Petersson metric in terms of limits of the second derivatives of length functions. In terms of geodesic currents Wolpert's description can be written as    
   $$L''_{\mu_{X}}(v, w) = <v ,w>_{g} \mbox{ for } v, w \in T_{X}(F(S))$$ (see \cite{Bon88}).
Therefore 
   $$||v||^{2}_{H} = ||v||^{2}_{g}$$ giving the result.
   \eproof

\end{document}